\documentclass[preprint]{elsarticle}
\usepackage[latin1]{inputenc} 
\usepackage{amsmath} 
\usepackage{amsfonts} 
\usepackage{amssymb} 
\usepackage{amsthm} 
\usepackage{graphicx}  
\usepackage{latexsym}    

\usepackage{amsmath}
\usepackage{graphicx}
\usepackage{amsfonts}

\newtheorem{corollary}{Corollary}[section]
\newtheorem{theorem}{Theorem}[section]
\newtheorem{lemma}{Lemma}[section]
\theoremstyle{definition}
\newtheorem{remark}{Remark}[section]
\newproof{pf}{Proof}

\newcommand{\bbc}{{\boldsymbol c}}
\newcommand{\bbw}{{\boldsymbol w}}
\def\rpn{\mathbb{R}_+^N}
\newcommand{\bw}{\boldsymbol w}
\newcommand{\bW}{\boldsymbol W}
\newcommand{\bc}{\boldsymbol c}
\newcommand{\bC}{\boldsymbol C}
\newcommand{\bx}{\boldsymbol x}
\newcommand{\bb}{\boldsymbol b}
\newcommand{\bq}{\boldsymbol q}
\newcommand{\bt}{\boldsymbol \tau}
\newcommand{\clC}{{\boc{C}}}

\usepackage{graphicx}   
\usepackage[english]{babel}   
   
\usepackage[usenames]{color}
 \newcommand{\bfj}{\boldsymbol{\hat f}}

\def\boc#1{\boldsymbol{{\cal #1}}}
\def\mD{\boc{D}}
\def\bh#1{\boldsymbol{\hat{#1}}}
\usepackage{xcolor}
\usepackage[normalem]{ulem}

\begin{document}

\begin{frontmatter}

\title{Implicit-Explicit WENO schemes for the equilibrium dispersive model of chromatography.}    
\author{Rosa Donat \footnote{Email: donat@uv.es}, Francisco Guerrero \footnote{Corresponding author. Email: guecor@uv.es Tel: +34963543232 Fax:+34963543922} and Pep Mulet \footnote{Email: mulet@uv.es}%
\thanks{Department of Mathematics. Universitat de Val\`encia. Av. Dr. Moliner, 50; 46100-Burjassot-Valencia. Spain}%
}

\begin{abstract}
Chromatographic processes can be modeled  by nonlinear,
convection-dominated partial differential equations, together with
nonlinear relations: the {\em adsorption isotherms}. In this paper we
consider the nonlinear {\em equilibrium dispersive} (ED) model with
adsorption isotherms of {\em Langmuir type}.  We show that very
efficient, fully conservative,  numerical schemes can be designed for
this mode by exploiting the relation between the {\em conserved
  variables} of the model and the physical concentrations of the
multi-component mixtures. We show that this relation   is one to one
and  
admits a smooth global inverse, which cannot be given explicitly but can be
easily computed by using a convenient root finder. These results
provide the necessary  ingredients to implement  fully conservative
numerical schemes  for the model considered. 

Implicit-Explicit (IMEX) techniques can be used in the convection-dominated
regime in order to increase the efficiency of the numerical scheme.
We propose  a second order IMEX scheme, combining an explicit  
Weighted-Essentially-non-Oscillatory discretization of the convective
fluxes with an implicit treatment of the diffusive term, in order to 
ilustrate  the numerical issues involved in
the application of IMEX techniques to this model. Through a series of numerical experiments, we show
that the  scheme provides
accurate numerical solutions which capture the sharp discontinuities
present in the chromatographic fronts, with the same  stability
restrictions as in the purely hyperbolic case.

\end{abstract}

\begin{keyword} 
Numerical methods \sep WENO schemes \sep Chromatography \sep Implicit-Explicit methods \sep Conservation laws. 
 \end{keyword}

\end{frontmatter} 


\section{Introduction} 
Chromatography is a powerful tool for the separation of 
 complex mixtures. In liquid batch
chromatography, a pulse of fluid mixture (the solute) is injected at 
one end of a long cylindrical column filled with a porous medium (the stationary
phase), followed by a continuous flow of liquid  (the mobile phase) 
along the column. The solute interacts with the porous medium and is
distributed between the liquid and solid phases, and the components of the mixture begin to separate according to the strength of their
interaction with the stationary phase.  For a sufficiently long column,
 band profiles of single component-concentration travel  along the
column and it is possible to collect pure fractions of components at
the outlet of the device. These tools are used for difficult separation
tasks when a high purity of the product is demanded, as it is often
the case in the pharmaceutical industry.

 It has been long recognized that
chromatographic processes can be modeled by considering non-linear,
convection-dominated partial differential equations \cite{Mazzotti2013,Guiochon}, coupled with
some algebraic relations between the concentrations of the components
of the mixture in the mobile and solid phases. Under reasonable
assumptions, such as negligible dispersion effects and transport
resistances, these equations become systems of first order non-linear
conservation laws. Understanding
the mathematical theory of these systems can enlighten many of the engineering
aspects  \cite{Mazzotti2013}, in particular the formation and evolution of shock
waves, which are an essential ingredient in the formation of band profiles of pure components. In addition, and since analytical solutions can  only be obtained in very
simple situations, it is important to develop  tools that are
able to perform accurate 
numerical simulations using these models. As observed in \cite{Seidel},
robust and reliable numerical techniques  can help  practitioners to
reduce the need for costly trial-and-error empirical experimentation.

In this paper we concentrate on the {\em equilibrium-dispersive} (ED
henceforth) model. This is an ideal
model  based on the following assumptions (see e.g. \cite{Javeed,Guiochon})

\begin{enumerate}
\item There is a permanent equilibrium between the solid and mobile
  phases at all positions in the column. 
\item The compressibility of the mobile phase is negligible and there is no interaction
  between the solvent (carrier) and the solid phase.
\item The porous medium in the column is homogeneous. Then, the
  adsorption process is uniform in time and axial direction. 
\item There are no radial concentration gradients in the column. 
\item Only axial dispersion causes band broadening. The column
  efficiency is characterized by an apparent axial dispersion
  coefficient $D_a$, related to the height of the column, $L$, the
  (constant) velocity of the mobile phase, $u$, and the number of
  theoretical plates $N_t$, see \cite{Guiochon}, through the following
  relation    
\begin{equation*}D_a=\frac{Lu}{2N_t},
\end{equation*}

\item   Any additional factor that could influence the adsorption
  behavior (such as the temperature) is neglected. \end{enumerate}

The  mass balance equation of the ED model involves the
concentrations of the $N$  components of
the mixture in the mobile phase, ${\boldsymbol c}=(c_1,\ldots,c_N)^T$,
and the solid phase,
${\boldsymbol q}=(q_1,\ldots,q_N)^T$, and takes the following form
\begin{equation}\label{eq:modelED}
\frac{\partial {\boldsymbol c}}{\partial t}+\frac{1-\epsilon}{\epsilon}\frac{\partial {\boldsymbol q}}{\partial t}+u\frac{\partial {\boldsymbol c}}{\partial z}=D_a \frac{\partial^2 {\boldsymbol c}}{\partial^2 z}
\end{equation}
where $\epsilon$ is the total porosity of the solid phase,  $t$ is the time and
$z$ the axial coordinate along the column, that is normalized to have
unit height, so that the top is at $z=0$ and the bottom at $z=1$.
 Under the assumptions listed above, the equilibrium relationship
  between the solid 
phase and liquid phase concentrations is given by the {\em adsorption
isotherm} $\bq=\bq(\bc)$, which  is usually a non-linear function
\cite{Guiochon}. Appropriate boundary conditions for this model are proposed
in \cite{Guiochon}: 
\begin{equation}\label{eq:bc2}
  \left. u{\boldsymbol c}-D_a \frac{\partial {\boldsymbol c}}{\partial z} \right|_{z=0} = u{\boldsymbol c}_{inj}(t), \qquad \left. \frac{\partial {\boldsymbol c}}{\partial z}\right|_{z=1}=0,
\end{equation}
for a known function $c_{inj}(t)$.

The form of the adsorption isotherm determines the mathematical
structure of the solutions to the ED model.  When dispersion is negligible, 
the model equations (\ref{eq:modelED}) and the algebraic relation
$\bq=\bq(\bc)$ form a system of nonlinear, first order partial differential
equations. The mathematical structure of the model for $N=1$,
i.e.~{\em single-component} chromatographic elution, has been
described in \cite{Mazzotti2013} for various types of adsorption
isotherms.  

In this paper we consider multi-component mixtures for which the
adsorption isotherms are of {\em  
  Langmuir type}, that is
\begin{equation} \label{eq:q-def}
q_i = \frac{a_i c_i}{1+\sum_{i=1}^N b_i c_i}, \qquad i=1,2,\dots,N.
\end{equation}
where $a_i>0$ are the Henry coefficients,  and the coefficients $b_i>0$
quantify the 
nonlinearity of the isotherm. For $N=1$ and $D_a=0$, the analysis of
the resulting hyperbolic conservation law carried out in
\cite{Mazzotti2013} shows that   the solutions are  
characterized by continuous or 
discontinuous composition fronts that propagate along the separation
unit. For $0<D_a<<1$, (\ref{eq:modelED}) becomes a parabolic,
convection dominated PDE
 whose solutions  may display very sharp fronts. The mathematical
 theory for the multi-component case seems to be much less developed.

 Numerical simulations 
involving the nonlinear 
system (\ref{eq:modelED}) require efficient numerical
techniques that can accurately describe discontinuous fronts. As
reported by various authors (see e.g. \cite{Javeed} and references
therein),  Finite 
element (FE) methods, normally used for diffusion dominated
problems, often lead to numerical  oscillations in convection dominated
problems whose solutions  display sharp gradients, and it is also well known
that spurious numerical
oscillations are  also observed when
classical finite difference schemes (FD) are used for such
problems.

 In \cite{Javeed}, the ED model
(\ref{eq:modelED})-\eqref{eq:q-def} is rewritten as 
\begin{equation} \label{eq:modelwc}
 \frac{\partial \bw}{\partial t} + \frac{\partial (u\bc)}{\partial z}= 
D_a \frac{\partial^2 \bc}{\partial^2 z}, \qquad
\bw=\bW(\bc)=\bc+\frac{1-\epsilon}{\epsilon} \bq(\bc)
\end{equation}
and the authors propose to use a {\em conservative} 
discretization of the convective terms, $\partial_z (u\bc)$, combined
with  a standard
  centered discretization of the parabolic terms,  in a finite volume
  (FV) framework.
This numerical technique  relies on {\em the
  understanding} 
that there is a one to one correspondence between the  
variables $\bw$ and the concentrations $\bc$, so that
  \eqref{eq:modelwc} becomes a system of conservation laws   when $D_a=0$. Then,
a {\em conservative
    discretization} of 
  the convective terms  guarantees mass conservation  for the
  {\em conserved    variables}, $\bw$, and, as a
  consequence, the  {\em 
    shock-capturing} property, i.e.~shocks (for $D_a=0$) or steep profiles (for
  $D_a>0$) in the numerical solution have the correct speed of
  propagation (the reader is referred to e.g. \cite{Levequeb} for a
  complete description of conservative schemes for systems of
  conservation laws).  

Numerical schemes that combine a conservative discretization of the
convective terms with a standard discretization of the parabolic terms
  have been successfully used to compute numerical approximations to
  the solution of convection-dominated second order PDEs and systems
  (see e.g. \cite{Donat1,BMV13,BBMRV15}).
   However, since  the
function  $\bC(\bw)$  cannot be explicitly determined  when $N>1$, 
  Javeed {\em et al.}
propose to update the values of the vector $\bbc$ by solving instead
the following linearized version of \eqref{eq:modelwc}
\begin{equation} \label{eq:dcdt}
(I+\frac{1-\epsilon}{\epsilon} \frac{\partial \bq}{\partial \bc})
\frac{\partial \bc}{\partial t} + 
u \frac{\partial \bc}{\partial z} = D_a \frac{\partial^2 \bc}{\partial
  z^2}
\end{equation}
 using an upwind, flux-limited, high resolution, conservative discretization
of the derivative of the convective flux $u \bc$. The approach
in \cite{Javeed} is  attractive because it
 incorporates 
 modern shock-capturing numerical techniques in the computation of
 the convective fluxes in (\ref{eq:dcdt}), leading to numerical
 solutions which are free of numerical oscillations. However, the need to update directly
 the vector $\bc$ forces the authors to abandon the {\em conservative
   formulation } \eqref{eq:modelwc}
 of the ED model and, as a consequence, we shall see that the resulting
scheme fails to be conservative,  leading to wrong speeds
in the propagating fronts.

One of the objectives of this paper is to show that there is a
globally well-defined,
one-to-one co\-rres\-pondence between  the vector
of concentrations $\bc$ and the {\em conserved variables} $\bw$,  so that \eqref{eq:modelwc} can
be rewritten as follows:  
\begin{align} \label{eq:modelwc1}
&\frac{\partial \bw}{\partial t} + \frac{\partial \boldsymbol{f}(\bw)}{\partial
  z}=D_a \frac{\partial^2 \bC(\bw)}{\partial^2 z}, \qquad
\boldsymbol{f}(\bw)=u \bC(\bw),\end{align}
with $\bC(\bw)$ a continuously 
differentiable function, satisfying
$\bC=\bW^{-1}$.  We show that, although there is no explicit expression for the
 function $\bC(\bw)$ for $N>1$, the value of   $\bC(\bw)$  for any
 $\bw\neq 0$, $w_i\geq 0$ can be
 determined by computing the only 
positive root of a well defined rational function. Hence,  the
necessary transfer of information required by a 
conservative scheme can be  carried out.

 In addition, we show that the structure of the Jacobian matrix
  $\bW'(\bc)$ can be computed in a rather straightforward manner by
  using the {\em secular equation}, as 
  in \cite{Donat3}. As a consequence, we show that all the
  eigenvalues of the  Jacobian matrix
  $\bC'(\bw)$ are strictly positive and pair-wise different, which
  allows us to prove  the  strict hyperbolicity of 
  the model when   $D_a=0$, and the well-posed character of
  \eqref{eq:modelwc}-\eqref{eq:q-def}  for $D_a>0$.

 These 
  results provide the theoretical background to 
  implement state-of-the-art,  fully
conservative numerical schemes 
for the ED model, with Langmuir-type adsorption isotherms. In this paper we
propose to use a second order Implicit-Explicit
 Runge-Kutta (IMEX-RK) scheme, that incorporates an off-the shelf
 Weighted-Non-Oscillatory (WENO) discretization of the convective flux
 terms. IMEX-RK schemes for convection-dominated parabolic PDEs
 combine the efficiency inherent to an implicit treatment of the
 second order derivatives (since the stability restrictions on the
 time step are the same as the CFL restriction that holds for
 $D_a=0$), with the robustness associated to a non-linear, high order
 conservative discretization of the convective derivative, which is
 treated in an  explicit fashion.

The paper is organized as follows: in section \ref{sec:Modeleq} we analyze the
mathematical structure of the ED model.  In particular, we prove that
the inverse function $\boldsymbol{C}(\boldsymbol{w})$ is globally well
defined and smooth in $\rpn:=\lbrace \boldsymbol{y}\in\mathbb R^N
  \colon y_i\geq 0\rbrace$.
Our analysis relies on the eigen-structure of the Jacobian matrix
$\bW'(\bc)$, which can be determined by
rewriting it as a rank-one perturbation of a diagonal matrix.  
In section \ref{sec:NumScheme}, 
we discuss the application of conservative schemes to the ED model
\eqref{eq:modelwc}-\eqref{eq:q-def}. We
show the effects of considering a non-conservative numerical scheme,
versus a fully conservative one, and discuss the convenience of using
implicit techniques when $D_a>0$. Section \ref{sec:imex-weno} describes  a simple second-order
IMEX-WENO scheme and discusses the various issues required for its
implementation in  numerical simulations of
chromatographic processes that fit the ED model
\eqref{eq:modelwc}-\eqref{eq:q-def}. In section \ref{sec:numex} we show some
numerical experiments to test the performance  of our WENO-IMEX-RK2
scheme. We close with some conclusions and perspectives for future work.

\section{The mathematical structure of the Equilibrium Dispersive Model}
\label{sec:Modeleq}

The ED model
\eqref{eq:modelED}-\eqref{eq:q-def} can be rewritten as
\eqref{eq:modelwc}, where  $\bw= \bW(\bc)=(W_i(\bc))_{i=1}^N$
with \begin{equation}\label{eq:definitionw}
  \begin{aligned}
    & \ W_i
    (\bc)= c_i \left( 1+\frac{\eta_i}{1+\sum_{j=1}^N b_j c_j} \right),
    \quad \eta_i=\frac{1-\epsilon}{\epsilon}a_i, \quad 1 \leq i\leq N
\end{aligned}
\end{equation}

Since $1+\sum_j b_j c_j \geq 1$, $\forall \bc \in \mathbb
R_+^N=[0,\infty)^N$, it follows that 
 the set $\mathbb R_+^N $ is in the interior
 of the 
 domain of the function $\boldsymbol{W}\colon \mathbb R^N \to
 \mathbb  R^N$. Moreover,  $\bW(\mathbb
R_+^N)\subseteq \mathbb R_+^N$.

We shall prove that  
$\boldsymbol{W}\colon \mathbb R_+^N \to  \mathbb  R_+^N$, 
is a continuously differentiable bijection.   The  {\em local } invertibility of
$\boldsymbol{W}$ will follow from the inverse function theorem, hence
we first analyze the Jacobian Matrix $\bW'(\bc)$. 
Since  this matrix can be written as a rank-one
perturbation of a diagonal matrix, the analysis of its structure can
be carried out via the {\em secular equation} (see also \cite{Donat3} and
references therein)

In what follows we shall assume that the components of the mixture are
ordered so that $0<a_1<a_2<\ldots<a_N$.

\begin{theorem} \label{th:jacobian}
For any  ${\boldsymbol  c}\in  (0,\infty)^N$,
  the Jacobian matrix ${\boldsymbol W}'({\boldsymbol c})$ is
  diagonalizable, with real, strictly positive, pairwise distinct, eigenvalues
  $\lambda_1,\dots,\lambda_N$ satisfying
  \begin{equation} \label{eq:Sdelam-roots}
 1< \lambda_1  <d_1 < \lambda_2  < \dots < d_{N-1} <
\lambda_N < d_N, \, \, d_i= 1+\eta_i(1+\sum_jb_j c_j)^{-1}.
\end{equation}

\end{theorem}
\begin{pf}
Consider the function $p:\rpn \rightarrow [1,\infty)$, 
$p({\boldsymbol c}):=1+\sum_jb_j c_j=1+\bb^T \bc$. 
 For any $\bc \in
\mathbb R_+^N$ we can write, with the aid
 of Kronecker's delta $\delta_{i,j}$
\begin{equation*}J_{ij}(\bc)=\partial_j W_i(\bc)=\delta_{i,j}
\left(1+\frac{\eta_i}{p(\bc)}\right)-\frac{\eta_i b_j c_i}{p(\bc)^2}. 
\end{equation*}
 Hence, if we define (dropping the explicit $\bc$ dependence for simplicity)
\begin{equation} \label{eq:D-def}
 d_i:= 1+\eta_i/p(\bc), \quad 
\tau_i:=-\frac{\eta_i c_i}{p(\bc)^2}, \quad i=1,\ldots, N
\end{equation}
we can write 
\begin{equation*}J=D+\bt \bb^T,   \qquad D:=\text{diag}(d_i)_{i=1}^N.
\end{equation*}

For any fixed $\bc \in (0,\infty)^N_+$, the
eigen-structure of such matrices can be easily determined
(see \cite{Donat3}) by computing the roots of the rational function
\begin{equation*}S_{\bc}(\lambda):=1+\bb^T(D-\lambda I)^{-1} \bt = 1+ \frac{1}{p(\bc)^2} \sum_{i=1}^N
\frac{\eta_i b_i c_i}{\lambda -d_i}.
\end{equation*}

Notice that  its poles, $d_i=1+\eta_i/p(\bc)>1$ satisfy $d_1<d_2<\cdots<d_N$ (since
$0<a_1<a_2<\ldots<a_N$). We can easily check that 
\[ S_{\bc}(\pm \infty)=1, \qquad \lim_{\lambda\to d_i^\pm}S_{\bc}(\lambda)=\pm
\infty, \qquad S_{\bc}(1)=\frac{1}{p(\bc)} >0. \]
  Hence, there must be at least $N$ roots of $S_{\bc}(\lambda)$,
$\lambda_i$, $1\leq i\leq N$,  satisfying
\begin{equation*} 
 1< \lambda_1 <d_1 < \lambda_2 < d_2 < \cdot\cdot\cdot < d_{N-1} <
\lambda_N < d_N. 
\end{equation*}
It is easy to see that these are the only roots of $S_{\bc}(\lambda)$, 
since  $S_{\bc}(\mu)=0$
implies $Q(\mu)=0$ for 
\[Q(\lambda)=S_{\bc}(\lambda)\Pi_{i=1}^N(\lambda-d_i),\] which is a
polynomial of degree $N$, with $N$ roots at most. Hence, the roots 
$\lambda_i$, $i=1,\ldots, N$ in \eqref{eq:Sdelam-roots} must be all
the roots of $Q(\lambda)$ and, as a consequence, all the roots of
$S_{\bc}(\lambda)$.

The above argument shows that $\forall \bc \in (0,\infty)^N$,  $S_{\bc}(\lambda)$ has $N$ different,
strictly positive roots. We claim that these roots   are, precisely, the 
eigenvalues of  the Jacobian matrix $J=\bW'(\bc)$.

To prove the claim above, we observe that $D-\lambda_i I$ is
invertible for any  of the roots of $S_{\bc}(\lambda)$,
hence, for $\tau$ in \eqref{eq:D-def}  we may  define ( $\bc\neq 0$
implies $\bt\neq 0$)  
\begin{equation*} \bx^{i}:= -(D-\lambda_i I)^{-1}\bt \neq 0.
\end{equation*}
and check that $J\bx^{i}=(D+\bt \bb^T)\bx^{i}=\lambda_i \bx^{i}$.
 Observe that 
\[0=S_{\bc}(\lambda_i)=1+\bb^T(D-\lambda_i I)^{-1}\bt= 1-\bb^T\bx^{i},
\]
which implies $\bb^T\bx^{i}=1$, 
hence
\[ J \bx^{i}= D\bx^{i} + \bt (\bb^T\bx^{i})= (D-\lambda_i I)\bx^{i}+\lambda_i \bx^{i} +\bt
= \lambda_i \bx^{i}. \]  

Thus, $J(\bc)={\boldsymbol   W}'({\boldsymbol c})$  has $N$ strictly positive,
pairwise distinct, eigenvalues and it is, therefore,  diagonalizable.
\end{pf}

\begin{remark}
Theorem \ref{th:jacobian} implies that
$\bW'(\bc)$ is  non
singular $\forall \bc \in \mathbb R^N_+$, hence  the inverse function theorem guarantees the existence
of a local inverse in a neighborhood of any $\bw=\bW(\bc)$, $\bc \in
\mathbb R^N_+$.
 We shall prove that, in fact,  there is a globally defined
inverse function $\bC: \mathbb R^N_+ \rightarrow \mathbb R^N_+$. For this, we
prove first the following result.
\end{remark}

\begin{lemma} \label{lemma:1}
For any fixed $\bw \in \rpn$, the rational function $R_{\bw}:
  \mathbb{R}  \rightarrow   \mathbb{R}$, 
\begin{equation} \label{eq:Rp-def}
R_{\bw}(y)=1-y+\sum_{i=1}^N \frac{y}{y+\eta_i} b_i w_i.
\end{equation} 
 has only one positive root, $\rho_0(\bw)$. In
addition,  $1 \leq \rho_0(\bw) \leq p(\bw)=1+ \bb^T \bw$.
\end{lemma}
\begin{pf}
  It is easy to see that $w_i>0$, $i=1,\dots,N$, may be assumed without loss of
    generality. 
 The  rational function $R_{\bw}(y)$ 
satisfies
\begin{equation*}
\lim_{y\to -\eta_i^\pm}R_{\bw}(y)=\mp \infty,  \qquad  
\lim_{y\to \pm \infty}R_{\bw}(y)=\mp \infty, \qquad R_{\bw}(0)=1
\end{equation*}
hence it has at least $N+1$ real roots,  $\rho_0, \ldots, \rho_N$, satisfying
\[   -\eta_N < \rho_N <  \cdot\cdot\cdot < -\eta_1 < \rho_1 < 0 <\rho_0\]
Note that for  $p(\bbw)=1+\bb^T \bbw$ we have
\[ R_{\bw}(p(\bw))= -\sum_i b_i w_i + \sum_i b_i w_i 
\frac{p(\bw)}{p(\bw) +  \eta_i} 
= -\sum_i b_i w_i \frac{\eta_i}{p(\bw) +
  \eta_i} <0 \]
while  $R_{\bw}(1)=\sum_i b_i w_i /(1+\eta_i) > 0$,
hence, we must have $1 \leq \rho_0 \leq 1+ \bb^T \bw$.
 
 Since any root of $R_{\bw}(y)$ is also a root of 
\[Q(y)=R_{\bw}(y)\, \Pi_{i=1}^N (y+\eta_i), \]
which is a polynomial of degree $N+1$, $\rho_0,\ldots,\rho_N$ must be all
the roots of $Q(y)$ and, hence, of $R_{\bw}(y)$.

\end{pf}

\begin{theorem}\label{th:winvertible}
The function ${\boldsymbol W}:
\rpn \rightarrow \rpn$ given by
  \eqref{eq:definitionw} is invertible. The  inverse function $\bC:=\bW^{-1}\colon \rpn \rightarrow \rpn$ is   continuously differentiable in $\rpn$ and is
defined by
\begin{equation} \label{eq:cw-def}
 C_i(\bw) := \frac{w_i}{1+{\eta_i}/{\rho_0(\bbw)}} 
\end{equation}
where $\rho_0(\bbw)$ is the only positive root of the rational function
$R_{\bw}(y)$ in \eqref{eq:Rp-def}.

\end{theorem}
\begin{pf}

Let $\bc=(c_1,\ldots,c_N)^T\in \rpn$ and
$p(\bbc):=1+\bb^T\bbc\geq1$. According to (\ref{eq:definitionw})\[ W_i(\bbc)=c_i \left(1+\frac{\eta_i}{p(\bbc)}\right) \qquad \equiv \qquad c_i=
\frac{W_i(\bc)}{1+\eta_i/p(\bbc)}, \]
hence we can write
\begin{align*}
  p(\bbc)-1=\sum_{j=1}^N b_j c_j = 
\sum_{j=1}^N \frac{ b_j  W_j(\bc)}{1+\eta_j/p(\bbc)}\end{align*}
or, equivalently,
\begin{align*}
0=1-p(\bbc)+p(\bbc)\sum_{j=1}^N \frac{ b_j
  W_j(\bc)}{p(\bbc)+\eta_j}=R_{\bW(\bc)} (p(\bc)),
\end{align*}
since $p(\bc)\geq1$,  $p(\bc)$ must be the only positive root of $R_{\bW(\bc)}(y)$ in
(\ref{eq:Rp-def}), i.e. $p(\bc)=\rho_0(\bW(\bc))$. 
Then, for any $\bc \in \rpn$, 
\[C_i(\bW(\bc))=\frac{W_i(\bc)}{1+{\eta_i}/{\rho_0(\bW(\bc))}}=
\frac{W_i(\bc)}{1+{\eta_i}/{p(\bc)}}=c_i,
\]
that is, $\bC\circ \bW= \text{id}_{\mathbb R_+^N}$.

Consider now $\bw\in\mathbb R_+^N$. Taking into account Lemma \ref{lemma:1}:
\begin{align*}
  &p(\bC(\bw))=1+\sum_{i=1}^{N} b_i C_i(\bw)=1+\sum_{i=1}^{N} b_i
  \frac{w_i}{1+\eta_i/\rho_0(\bw)} \\
  &=1+\rho_0(\bw)\sum_{i=1}^{N} b_i
  \frac{w_i}{\eta_i+\rho_0(\bw)} =
  R_{\bw}(\rho_0(\bw))+\rho_0(\bw)=\rho_0(\bw).
\end{align*}
With this result, it easily follows that $\bW\circ \bC=
\text{id}_{\mathbb R_+^N}$  since
\begin{align*}
  W_i(\bC(\bw))=C_i(\bw)\left(1+\frac{\eta_i}{p(\bC(\bw))}\right)=
  w_i\left(1+\frac{\eta_i}{\rho_0(\bw)}\right)^{-1}\left(1+\frac{\eta_i}{p(\bC(\bw))}\right)=w_i.
\end{align*}  

It follows from the inverse function theorem that the function $\bC$
is continuously differentiable in  $\rpn$.
\end{pf}

\begin{corollary}\label{cor:1}
  The ED model \eqref{eq:modelwc1} is well
  posed. For $D_a=0$ the system of conservation laws  is strictly
  hyperbolic, and for any $\bw$ such that $w_i>0$, all the
  eigenvalues $\mu_j$ of the Jacobian   matrix   $u C'(\bw)$ are positive,
  pairwise distinct, and bounded above by $u$. With the
  notation in Theorem \ref{th:jacobian}, the
  $\mu_j=u/\lambda_j$ satisfy
  \begin{equation*}  u > \mu_1 >u/d_1 > \mu_2  > \dots > u/d_{N-1} >
\mu_N > u/d_N > 0.
\end{equation*}

For $D_a>0$, the system
  is parabolic in the sense of Petrovskii
  (cf. ~\cite{EidelmanZhiharatshu98}), i.e., the eigenvalues of the
  matrix $D_a\bC'(w)$ are bounded below by some positive constant for
  any $w\in\rpn$.  
\end{corollary}
\begin{pf}
 Theorem \ref{th:winvertible} allows us to write  system \eqref{eq:modelwc} as
\begin{equation} \label{eq:modelwf}
\frac{\partial {\boldsymbol w}}{\partial t}+\frac{\partial (u
  {\boldsymbol C}(\boldsymbol w))}{\partial z}=D_a
\frac{\partial}{\partial z}\left[\boldsymbol C'(\boldsymbol
  w)\frac{\partial {\boldsymbol w}}{\partial z}\right] 
\end{equation}

According to Theorem \ref{th:jacobian}, all the eigenvalues of the
Jacobian matrix $\boldsymbol{W}'(\boldsymbol{c})$ are positive and
greater than 1, i.e. $\lambda_i >1$ for $i=1,2,...,N$. This implies
that the eigenvalues of the Jacobian matrix
$\boldsymbol{C}'(\boldsymbol{w})$ for the inverse function are: 
\[ \frac{1}{\lambda_i} < 1, \quad \mbox{for}\quad i=1,2,...,N. \]
Since in the ED model the flux is given by:
\[ \boldsymbol{f(w)} = u \boldsymbol{C}(\boldsymbol{w}), \]
all the characteristic speeds of the Jacobian matrix of the physical
flux satisfy the bounds in the statement from
\eqref{eq:Sdelam-roots}.

 For $D_a>0$, the matrix $D_a\boldsymbol
 C'(\boldsymbol w)$ appearing in the  diffusion term,
 has eigenvalues $D_a\lambda_i(\bC'(\bw))^{-1}$ from Theorem
 \ref{th:jacobian} and Theorem \ref{th:winvertible}. From
\eqref{eq:Sdelam-roots} these eigenvalues are bounded below by
\begin{align*}
  \frac{D_a}{ d_N}=\frac{D_a}{1+\frac{\eta_N}{1+\sum_jb_j c_j}}\geq
  \frac{D_a}{1+\eta_N} >0
\end{align*}
and therefore \eqref{eq:modelwf} is a
parabolic system  in the sense of Petrovskii.

\end{pf}

\begin{remark}\label{remark:1}
 The  function ${\boldsymbol C}({\boldsymbol  w})$
cannot be written 
explicitly for $N>1$. However, for all practical purposes, the  computation of
$\bC(\bbw)$ in \eqref{eq:cw-def} only requires to  find  $\rho_0(\bw)$,
   the only positive root of
the rational function \eqref{eq:Rp-def}.
Since $1 \leq \rho_0(\bw)\leq 1+\bb^T \bw$,   this root can easily be   found using a convenient  root finder. 
\end{remark}

\section{Conservative Numerical Schemes}\label{sec:NumScheme}

As mentioned in previous sections, discontinuous fronts (for $D_a=0$)
or sharp profiles (for $0<D_a<<1$) are to be  expected when computing
numerical approximations to \eqref{eq:modelwc}. 
 It is well known that such solutions  are notoriously hard 
 to compute accurately with 
classical numerical schemes, which tend to produce moving fronts that
display oscillatory, Gibbs-like, behavior.  Conservative,
shock-capturing, schemes respect a fundamental property: the conservation 
in time of the 'total mass' of  the {\em conserved variables} ($\int
\bw(x,t) dx$
in the ED model). This fact (via the Lax-Wendroff
theorem, see \cite{Levequeb} for further details), leads to  
numerical solutions which have  the correct behavior in terms of the
velocity  of shocks.  Combined
with appropriate 
discretization of the diffusive terms, a high resolution conservative
discretization of the convective terms can be used to 
'capture' the sharp fronts that appear in convection dominated
systems, avoiding the numerical oscillations observed in more
conventional schemes.

In \cite{Javeed}, the first attempt at
using  modern shock capturing techniques for the model equation
\eqref{eq:modelED} is carried out.  The starting point in \cite{Javeed} is the following 
semi-discrete scheme for the evolution of the cell-averages of $\bw$
in  a
uniform mesh with grid points $z_j=(j-\frac{1}{2})\Delta z$,
$j=1,\dots,m$, where $\Delta z=L/m$, $L$ being the length of the
column ($z=0$ and $z=L$ represent the beginning and the end of
the column, respectively) 
\begin{equation}\label{eq:cons_form}
\frac{d {\boldsymbol w}_j (t)}{d t}=-\frac{1}{\Delta z} 
 \left(
{\boldsymbol{ \hat f}}_{j+1/2}-{\boldsymbol{\hat f}}_{j-1/2} 
\right)+ \frac{D_a}{\Delta z} \left ( 
\left( \frac{\partial c}{\partial z}\right )_{i+1/2}-
\left( \frac{\partial c}{\partial z}\right )_{i-1/2} \right ).
\end{equation}
Here 
$\bw_j(t)\approx \frac{1}{\Delta
  z}\int_{z_{j-1/2}}^{z_{j+1/2}}\bw(z,t)dz$  
and $\hat{\boldsymbol{f}}_{j+1/2}= \hat{\boldsymbol{f}}({\boldsymbol
  w}_{j-p},\dots,{\boldsymbol w}_{j+q})$, for some function
$\hat{\boldsymbol{f}}$ of $p+q+1$ arguments, 
is the numerical flux at  the cell-boundary $z_{j+1/2}$. As observed
in \cite{Javeed}, 
the simplest choice is the first order {\em upwind} numerical 
flux, which for the ED model becomes $\bfj_{j+1/2} = u\bc_j $ (notice that our analysis confirms that  all propagation speeds
are positive, hence this is indeed the {\em upwind} choice). 

 Considering the first order upwind numerical fluxes,  centered differences for the parabolic term and 
  the Forward Euler method  for the time evolution, we get the
 following  conservative, fully discrete scheme
\begin{equation}\label{eq:fullydiscreteFE}
{\boldsymbol w}_j^{n+1} ={\boldsymbol w}_j^{n}-\frac{\Delta t}{\Delta
  z} \left(
u \bc_j - u \bc_{j-1} \right)
+\frac{D_a \Delta t}{\Delta z^2} \left(\bc_{j+1} -2 \bc_j + \bc_{j-1}\right).
\end{equation}
which is  first order in space and time. Higher order conservative  schemes may be
obtained from  \eqref{eq:cons_form} by using appropriate ODE solvers
for the time evolution  and  high-resolution numerical flux functions
in the discretization of the convective derivative.

The  lack of an explicit expression for the relation $\bbc=\bC(\bbw)$
prevented the authors in \cite{Javeed} from using directly the
conservative formulation of the numerical schemes derived from
\eqref{eq:cons_form}. As an alternative,  they consider  the following 
(non-conservative) linearization of  \eqref{eq:modelwc} 
\begin{equation*} 
(I+\frac{1-\epsilon}{\epsilon} \frac{\partial \bq}{\partial \bc})
\frac{\partial \bc}{\partial t} + 
u \frac{\partial \bc}{\partial z} = D_a \frac{\partial^2 \bc}{\partial
  z^2},
\end{equation*}
which allows the authors  to construct a numerical
scheme which updates directly 
the values of $\bc$.  A second-order (in space and time) scheme is
implemented in \cite{Javeed} by using flux-limiting techniques to construct a
high-resolution  numerical flux function and an appropriate ode solver.

The second order scheme proposed in \cite{Javeed} produces the  sharp,
non-oscillatory numerical profiles that characterize a
state-of-the-art high resolution scheme, however,  it  
fails to be conservative and the fundamental property of conservation of
the total mass no longer holds. As a consequence,    the shock
fronts (or the nearly discontinuous profiles, when $D_a<<1$) do not
move at the correct speed. 

On the other hand, according to the results of
the previous section, finding
 $\bc=\bC(\bw)$, only requires the computation of
 the  positive 
 root of a given rational function, so that the fully conservative
 version of these schemes can be easily implemented. 
 To illustrate the different behavior of the conservative versus the
 non-conservative numerical schemes, we consider the
ED model for single-component chromatographic elution:

\subsection{Single-component chromatographic elution.}\label{subsec:single}

 The model becomes (we do not need the vector notation) 
\begin{equation*}\frac{\partial  w}{\partial t}+\frac{\partial (u  c(w))}{\partial z}=D_a \frac{\partial^2  c(w)}{\partial z^2}.
\end{equation*}
For the Langmuir's isotherm:
\begin{equation*}w(c)=c \left(1+\frac{\eta }{1+bc}\right), \quad \eta=
\frac{1-\epsilon}{\epsilon} a
\end{equation*}
we can readily find
the  analytical expression for $c(w)$:
\begin{equation*}c(w)=\frac{1}{2b} \left(\sqrt{(1+\eta -bw)^2+4bw}-(1+\eta  -bw)\right).
\end{equation*}
which, as observed in Theorem \ref{th:winvertible},  corresponds to 
\[ c(w)=\frac{w}{1+\frac{\eta }{\rho_0(w)}}. \]
with $\rho_0(w)$ the unique  positive solution of the rational function
\[ R_w(y)=1-y+\frac{y bw}{y+\eta }=0. \]

For $D_a=0$, it is easy to see that  the quantity
\[ W(t)=\int_0^L w(c(z,t)) dz=\int_0^L \left( c(z,t)+\frac{1-\epsilon}{\epsilon} q(z,t)\right) dz\]
must be conserved at the continuous level during the  time evolution (after the moment when
all the solute has been injected and until it begins to leave the
column).
At the discrete level, we measure the conservation of total mass in the numerical solution  by computing the corresponding discretized magnitude:
\begin{equation} \label{eq:consWn}
 W^n:=\sum_j w_j^n\Delta z =\sum_j w(c_j^n)\Delta z = \sum_j \Delta z \left(
  c_j^n+\frac{1-\epsilon}{\epsilon} q_j^n\right).
\end{equation}
 Conservative schemes for hyperbolic conservation laws and systems
 maintain this quantity for all time.
 The following  test problem
shows that this
property no longer holds for the numerical schemes proposed in \cite{Javeed}.

We consider the following set of parameters $D_a=0$, $a=1, b=1$, $u=1$,
$\epsilon=0.5$ and assume that the component is injected between
$t=0$ and $t=0.2$ with $c=1$ at $z=0$. 

\begin{figure}[ht]   
\begin{center} 
\includegraphics[scale=0.35]{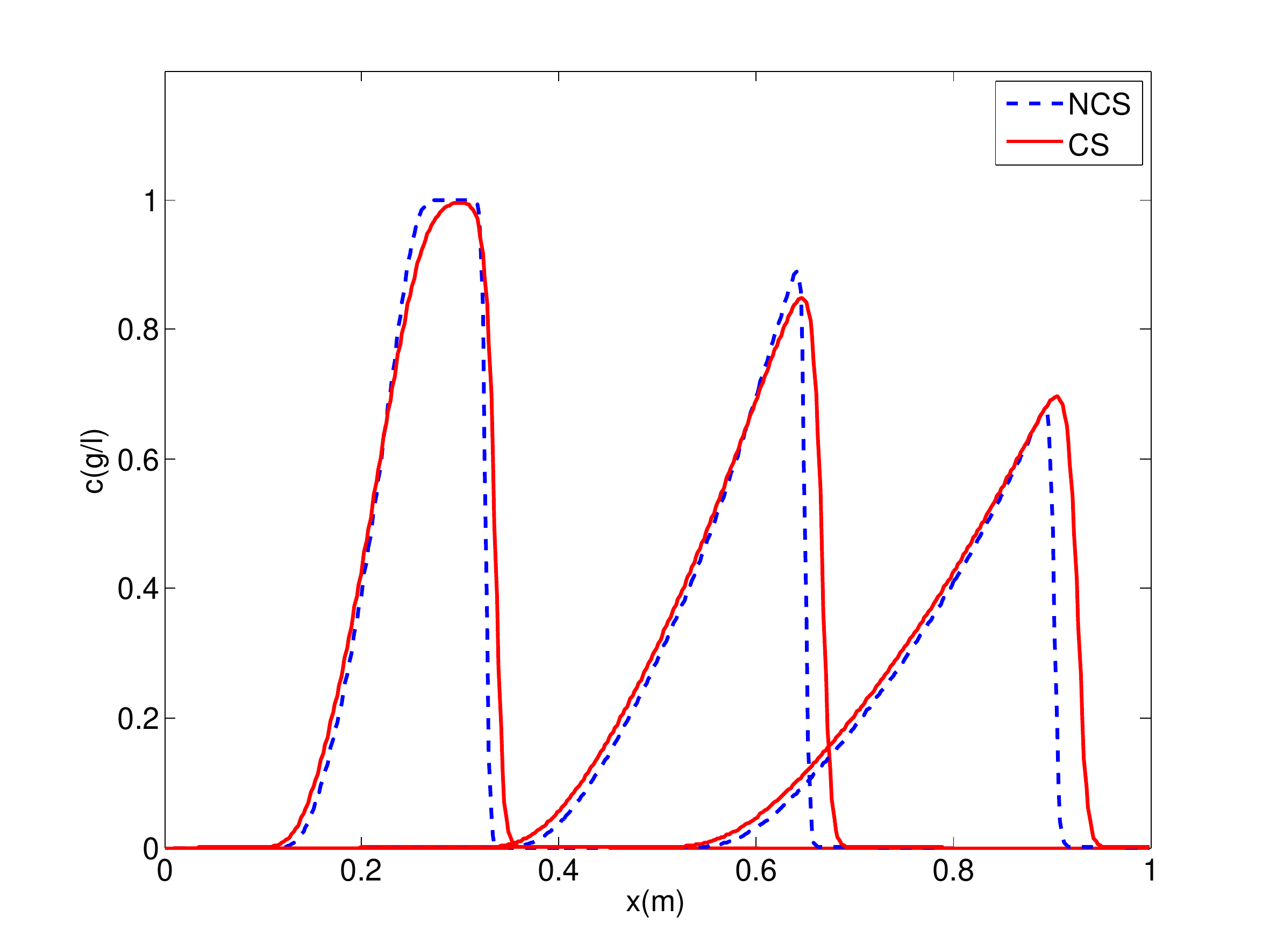}
\end{center}   
\caption{Single-component elution, $D_a=0$. Numerical solutions with conservative and non-conservative
  first order schemes. $m=500, k/h=0.9$, $T=0.5$, $T=1.0$ and $T=1.4$. }  
\label{comparison1}  
\end{figure}

\begin{figure}[ht]   
\begin{center} 
\includegraphics[scale=0.25]{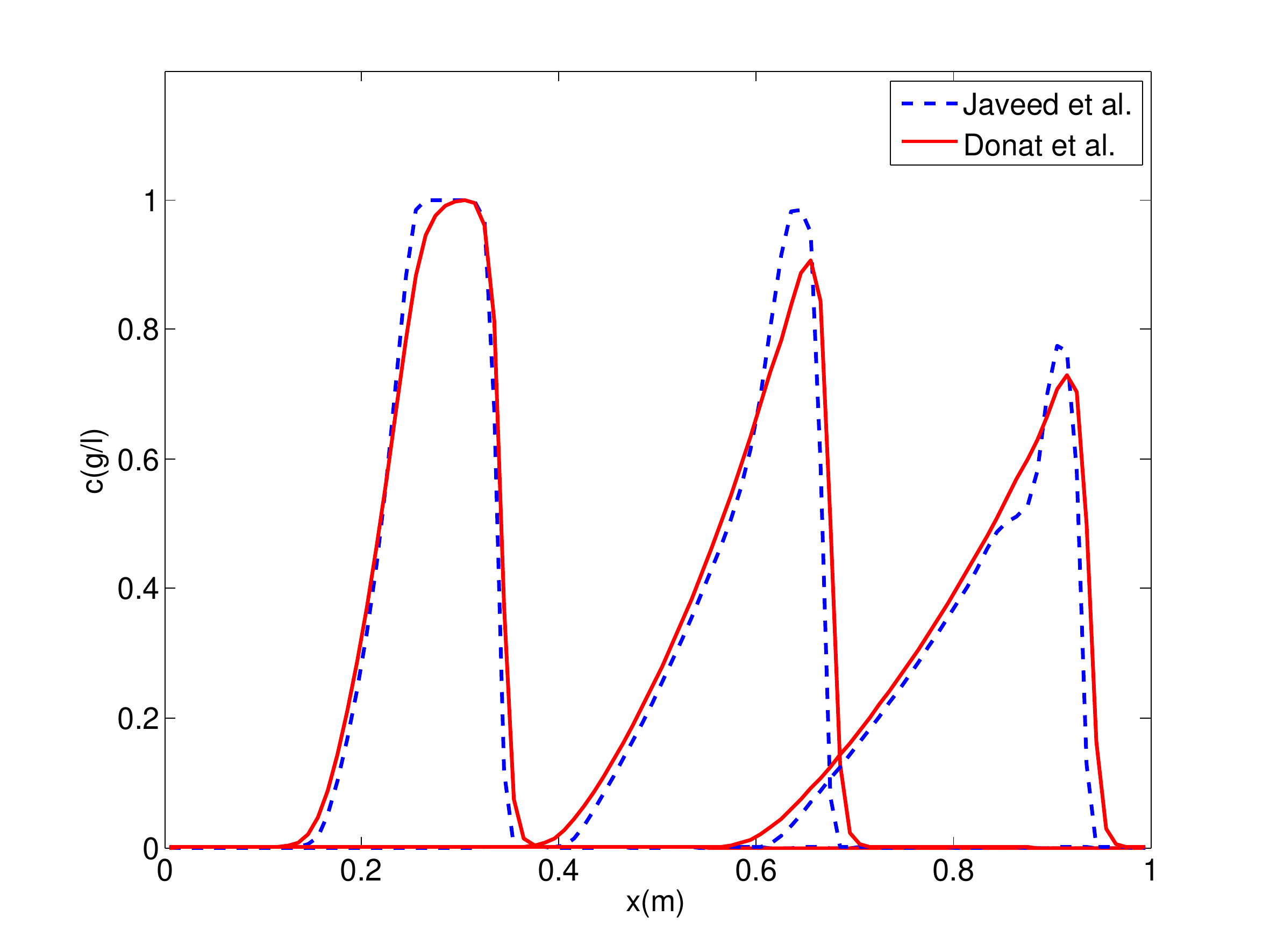}
\includegraphics[scale=0.25]{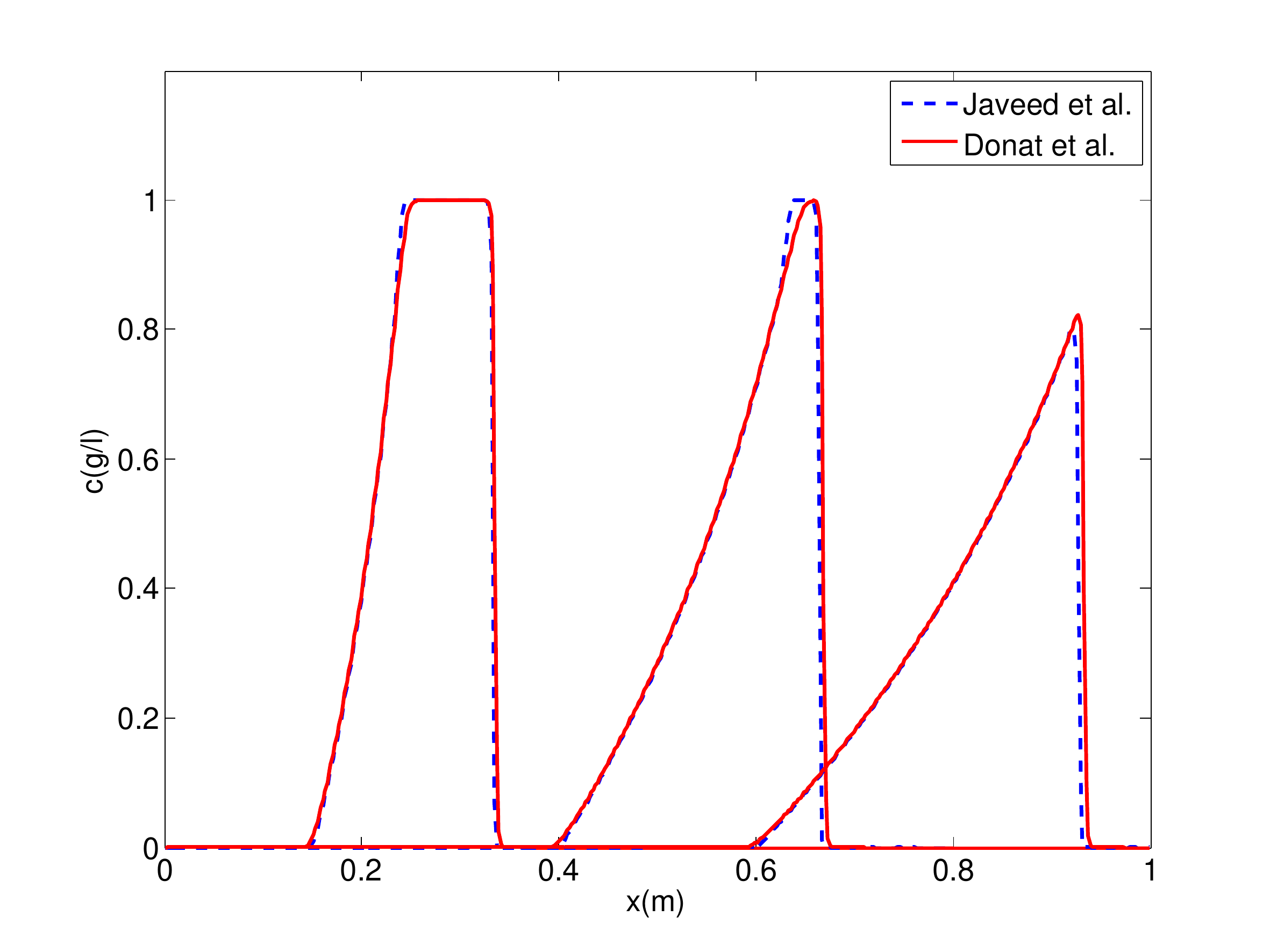}
\end{center}   
\caption{Single-component elution, $D_a=0$. Second order CS and NCS.  $T=0.5, 1.0, 1.4$. Left:
  $m=100$. Right:  $m=500$. $k/h=0.9$  }  
\label{comparison2}  
\end{figure}

Table \ref{tablacons} displays the
value of $W^n$ in \eqref{eq:consWn}, corresponding to three different 
times, obtained by following the non-conservative strategy (NCS) proposed in
\cite{Javeed} and a conservative scheme (CS) of the same order.
Figures \ref{comparison1} and \ref{comparison2} show snapshots of the numerical solution at the times
shown in the table, where the side  effect of the lack of conservation
can be clearly observed: when using  the non-conservative scheme,
the speed of the shock front is  slower than 
expected. 
\begin{table}[htbp]
\begin{center}
\begin{tabular}{|l|c|c|c|c||c|c|c|c|}
\hline
  & \multicolumn{2}{c|}{NCS} & \multicolumn{2}{c||}{CS}
& \multicolumn{2}{c|}{NCS} & \multicolumn{2}{c|}{CS} \\
\hline
 m & 100& 500 & 100 & 500 & 100& 500 & 100 & 500\\ \hline \hline
$T=0.5$ & 0.184 & 0.186 & 0.207  & 0.202  & 0.203 & 0.198 & 0.207 & 0.202  \\ \hline
$T=1.0$ & 0.171 & 0.173 & 0.207 & 0.202  & 0.200 & 0.194 & 0.207 & 0.202 \\ \hline
$T=1.4$ & 0.166 & 0.166 & 0.207 & 0.202  & 0.197 & 0.192 & 0.207
&0.202 
\\ \hline \hline
  & \multicolumn{4}{c|}{First order, $\Delta t/\Delta z=0.9$} & \multicolumn{4}{c|}{Second
    Order, $\Delta t/\Delta z=0.9$} \\ \hline
\end{tabular}
\caption{$D_a=0$. Conservation of  total mass for the single-component elution
  test.}
\label{tablacons}
\end{center}
\end{table}

The experiments carried out for the single elution experiment point
out that the  use of  small mesh 
parameters in the non-conservative formulation helps to mask the effect of
the loss of total mass. However, reducing the spatial step-size is not an
efficient alternative, due to the stability restrictions imposed by the
explicit treatment of the parabolic terms. A
 linearized Von-Neumann stability analysis
for the first order scheme \eqref{eq:fullydiscreteFE}  leads to
the following stability constraint (see \cite{Donat1} for details)
\begin{equation*}\frac{u\Delta t}{\Delta z}\max_w \varrho
(\boldsymbol{C'(w)})+\frac{2{\Delta t}D_a}{(\Delta z)^2} \max_w \varrho
(\boldsymbol{C'(w)})\le C_0 \le 1,
\end{equation*}
with $\varrho (\boldsymbol{C'(w)})$ the spectral radius of
$\bC'(\bw)$. Since
\[  \max_w \varrho (\boldsymbol{C'(w)}) < 1, \]
  numerical stability is obtained provided that
 \begin{equation}\label{conditionstability}
 \frac{\Delta t}{\Delta z} \le 
\frac{C_0}{u+\frac{2D_a}{\Delta z}} \quad \equiv \quad  \Delta t\leq \frac{C_0}{u\Delta z+2D_a}\Delta z^2.
\end{equation}
which implies a limitation on the time step
  $\Delta t \leq C_1 \Delta z^2$ when $u\Delta z < 2 D_a$.
It is shown in \cite{Donat1} that similar stability restrictions apply to
higher order explicit schemes for convection-diffusion equations and
systems.

\begin{figure}[ht]   
\begin{center} 
\includegraphics[scale=0.35]{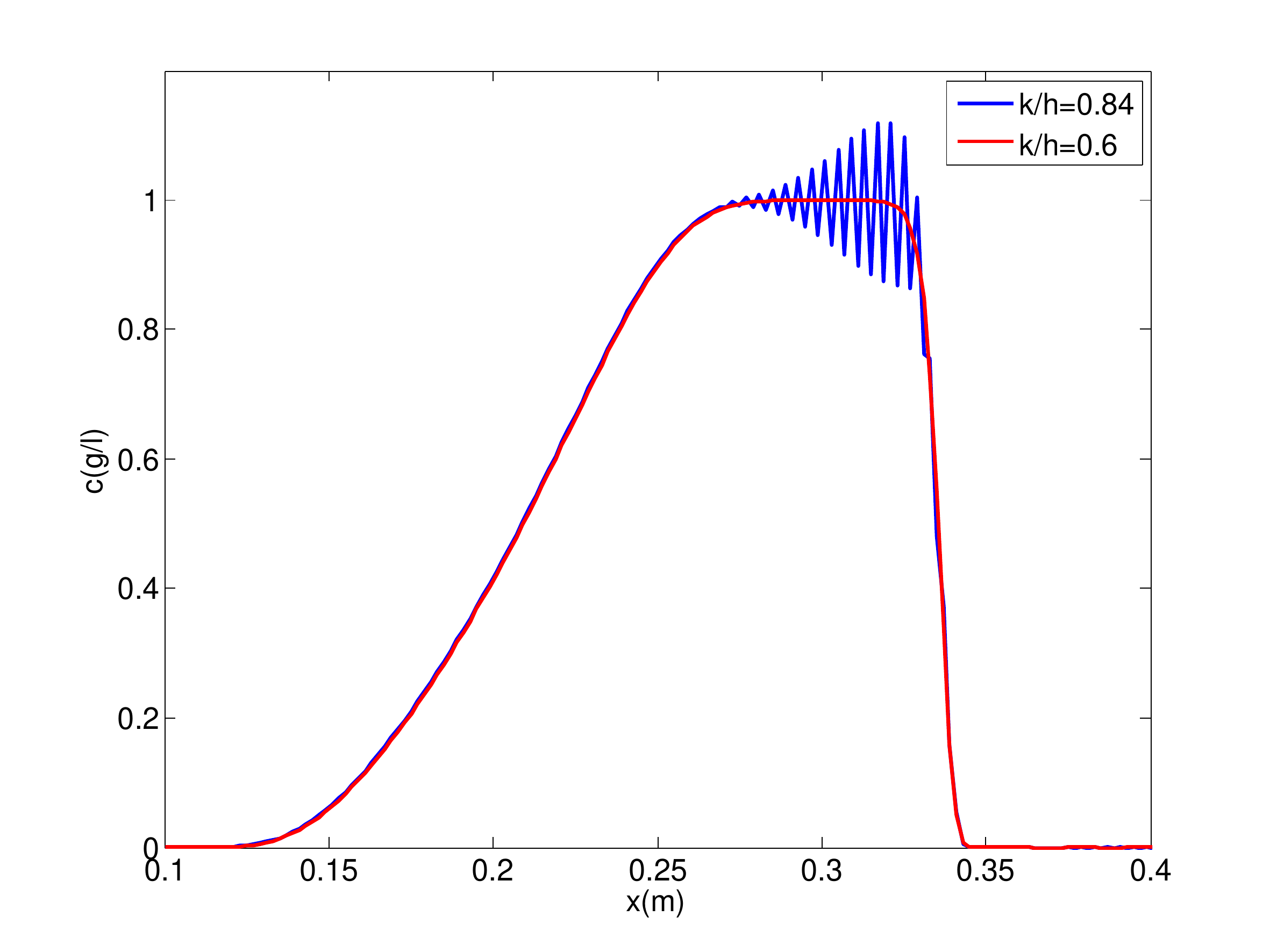}
\end{center}   
\caption{Single-component elution, $D_a=0.0005$. Numerical solutions with second
  order explicit CS . $m=500$, $T=0.5$ }  
\label{fig:stabCS1}  
\end{figure}

In Figure \ref{fig:stabCS1}  we
display the numerical solution at $T=0.5$ obtained with the second
order CS applied to the single-elution ED model
with  $D_a=0.0005$. The simulation corresponds to  $m=500$, for which
the stability bound in \eqref{conditionstability} becomes (taking
$C_0=1$) 0.6666. As observed in the figure, values of the ratio
$\Delta t/\Delta z$ above the stability bound produce
numerical oscillations that grow in time. The plot in Figure
\ref{fig:stabCS1} should be compared with the plots shown in Figure
\ref{comparison2}, that shows the simulation corresponding to $D_a=0$, for which
a ratio of $\Delta t/\Delta x=0.9$ is appropriate.

 The stability restrictions
that result from  the explicit treatment of the pa\-ra\-bo\-lic terms impose 
time steps that  can be   much smaller 
than the those required in the purely hyperbolic case. 
In order to avoid these strict stability restrictions  we
may turn to Implicit-Explicit strategies,  in which the parabolic term is handled
implicitly, while the convective term is discretized using any
convenient high-resolution shock capturing scheme. In  \cite{Donat1}
 we proved that the stability restriction for the first order upwind IMEX
 scheme, applied to  a convection-diffusion system similar to \eqref{eq:modelwf}
  is the same as for the purely hyperbolic case.

 Applying IMEX techniques to a convection-diffusion system requires to
 solve a system of equations at each time step, but very often these
  systems have a sparse structure and the cost of its solution is
  offset by the increased time step allowed by their stability
  constraints (see e.g. \cite{Donat1}). In the following section, we
  examine the specific issues that arise in the application of
  IMEX-WENO schemes for the ED model.


\section{IMEX-WENO schemes}\label{sec:imex-weno}
As observed in the previous section, Implicit-Explicit strategies
represent  an interesting option for the  numerical simulation of
chromatographic processes.  IMEX-WENO schemes have been
used in fairly similar scenarios \cite{Donat1,Guerrero} and extensive
numerical testing shows that they provide an efficient and robust
tool for the numerical simulations of convection-dominated parabolic
PDEs.  In this  
section we describe a simple  second order IMEX-WENO
strategy in order to illustrate the numerical issues involved in IMEX numerical
schemes for convection-diffusion systems, and in particular for the ED model.

For the application of this technique to the ED model, it is
appropriate to rewrite \eqref{eq:modelwc1} as follows
\begin{equation}
\label{eq:modelwc11}
\frac{\partial \bw}{\partial t} + \frac{\partial}{\partial
  z}\left(\boldsymbol{f}(\bw)- \boldsymbol{g}\left(\bw,
  \frac{\partial \bw}{\partial z}\right)\right)=0, \qquad
\boldsymbol{g}(\bw,
  \frac{\partial \bw}{\partial z})=D_a \frac{\partial \bC(\bw)}{\partial z},
\end{equation}

Considering the
same computational setting as before, i.e. $z_j=(j-1/2)\Delta z$,
the starting point in  a
high order 
WENO scheme applied to the ED model is the following   semi-discrete scheme
\begin{equation} \label{eq:weno-sd}
\boldsymbol{w}'(t)=\boldsymbol{{\cal{L}}}(\bw(t))+
\boldsymbol{{\mD}}(\bw(t))
\end{equation}
where $\bw(t)$ is an $m\times N$ matrix whose $j$-th column,
  $\bw_j(t)$, is an approximation of $ w(z_j, 
  t)\in\mathbb R^N$, $j=1,\dots,m$,
 $\boldsymbol{{\cal{L}}}$ represents the spatial discretization of the
convective term $-\frac{\partial}{\partial
  z} \boldsymbol{f}(\bw)$ and $\boldsymbol{{\mD}}$ the spatial
discretization of the diffusion term $\frac{\partial }{\partial
  z}\boldsymbol{g}\left(\bw,
  \frac{\partial \bw}{\partial z}\right)$ in
\eqref{eq:modelwc11}.

The schemes we propose are Finite Difference
schemes, hence, they compute numerical 
approximations to the point-values of the conserved variables,
$\bw_j(t_n)\approx \bw(x_j,t_n)$, $t_n = n \Delta t$ and are
characterized by a conservative discretization of the convective and
diffusive terms
of the form, (dropping the dependencies for simplicity)
\begin{align*}
\boldsymbol{{\cal L}}_j &= -\frac{1}{\Delta z} \left(
  \boldsymbol{\hat{{f}}}_{j+1/2}-\boldsymbol{\hat{{f}}}_{j-1/2}\right)\\
\boldsymbol{{\cal D}}_j &= \frac{1}{\Delta z} \left(
  \boldsymbol{\hat{{g}}}_{j+1/2}-\boldsymbol{\hat{{g}}}_{j-1/2}\right),
\end{align*}
using convective and diffusive numerical fluxes,
  $\boldsymbol{\hat{{f}}}_{j+1/2}$, $\boldsymbol{\hat{{g}}}_{j+1/2}$,
  respectively, that approximate the respective fluxes at the
  corresponding cell interface $z_{j+1/2}=z_j+\Delta z/2$.

The convective numerical flux
\begin{equation}\label{eq:203}
  \bfj_{j+1/2}=\bfj(\bw_{j-p},\cdots,\bw_{j+p+1})
\end{equation}  
is computed by finite
difference WENO schemes of order $2p+1$
(nowadays an almost black-box routine)
which entail  applying WENO reconstructions to split convective fluxes
\begin{equation*}
  f^{\pm}(w)=\frac{1}{2}\left(f(w)\pm \alpha_{j+\frac{1}{2}} w\right).
\end{equation*}  

We refer the interested reader to \cite{Shu2009} for more details, and simply mention here that
its blind application requires to specify a numerical viscosity
$\alpha_{j+\frac{1}{2}}$ at each
interface, which must be a local upper-bound of the size of all  the
characteristic speeds of the Jacobian matrix of the physical flux. For
the ED model, we know that all characteristic speeds are bounded by
$u$ (see Corollary \ref{cor:1}), hence we may take $u$ as the numerical viscosity at each
cell interface. In our numerical experiments we shall use the  WENO5 
numerical flux (in a  
component-wise fashion \cite{Shu2009}). We recall
  Remark \ref{remark:1} for the computation of the convective fluxes.

The diffusive numerical fluxes are computed by second order
  centered finite differences
  \begin{align*}    \bh{g}_{j+1/2}:=\bh{g}(\bw)_{j+1/2}=D_a\frac{\bC(\bw_{j+1})-\bC(\bw_{j})}{\Delta
    z}.
  \end{align*}

The discretization of the boundary conditions \eqref{eq:bc2} at
  $z=0=z_{1/2}$ is
given by prescribing the sum of convective and diffusive numerical
fluxes as follows:
\begin{align*}  \bh{f}_{1/2}-\bh{g}_{1/2}=u \bc_{inj}(t)=\left.f(\bw)-g\left(\bw, \frac{\partial \bw}{\partial
    z}\right)\right|_{z=0}.
\end{align*}
The term $\mathcal{L}_{1}+\mathcal{D}_{1}$ is modified accordingly:
\begin{align*}
\boldsymbol{{\cal L}}_1+\boldsymbol{{\cal D}}_1 &= \frac{1}{\Delta z} \left(
  \big(-\bh{f}_{3/2}+\bh{g}_{3/2}\big)-\big(-\bh{f}_{1/2}+\bh{g}_{1/2}\big)\right)\\
&= \frac{1}{\Delta z} \left(
  -\bh{f}_{3/2}+\bh{g}_{3/2}+u\bc_{inj}(t)\right).
\end{align*}

The discretization of the boundary conditions \eqref{eq:bc2} at
  $z=1=z_{m+1/2}$ consists in taking  $\bh{g}_{m+1/2}=0=\left.g\left(\frac{\partial \bw}{\partial
    z}\right)\right|_{z=1}$, so the term
$\mathcal{L}_{m}+\mathcal{D}_{m}$ is modified accordingly:
\begin{align*}
\boldsymbol{{\cal L}}_{m}+\boldsymbol{{\cal D}}_{m} &= \frac{1}{\Delta z} \left(
  \big(-\bh{f}_{m+1/2}+\bh{g}_{m+1/2}\big)-\big(-\bh{f}_{m-1/2}+\bh{g}_{m-1/2}\big)\right)\\
&= \frac{1}{\Delta z} \left(
  \bh{f}_{m-1/2}
  -\bh{f}_{m+1/2}-\bh{g}_{m-1/2}\right).
\end{align*}

With all the previous comments,
\begin{align*}
  \mD(\bw)&=\boldsymbol{C}^*(\bw)\mathcal{A},\quad
  \boldsymbol{C}^*(\bw)_{i,j}=C_i(\bw_j),
\end{align*}
where $\bw_j$ is the $j$-th column of the $N\times m$ matrix $\bw$ and
$\mathcal{A}$ is  the tridiagonal
$m\times m$ matrix  given by
\begin{align}\label{eq:207}
\mathcal{A}_{i,j}&=\begin{cases} -\mu & i=j=1,m\\
    -2\mu & i=j\neq 1,m\\
    \mu & |i-j|=1
  \end{cases},\quad \mu=D_a/\Delta z^2.
\end{align}

Furthermore, the needed convective fluxes as in \eqref{eq:203}, for
$j=1,\dots,m$ such that $k=j-p<1$ or $k=j+p+1 > m$, require values at the
corresponding ghost cells $x_k$, whose indices are, for the first case,
$k=-p+1,\cdots,0$ 
and for the second one, $k=m+1,\dots,m+p$. We  obtain the values at
those ghost cells by using  extrapolation with a linear polynomial
that satisfies the boundary condition and that interpolates the data 
for the internal point which is symmetric with respect to the
boundary. For $z=0$ and $k=1-j$, $j=1,\dots,p$, taking into account
\eqref{eq:bc2}, this extrapolation yields the value
\begin{align*}
  \bc_{1-j}=\frac{(-D_a/u-(j-1/2)h)\bc_j+2(j-1/2)h \bc_{inj}
  }{D_a/u+(j-1/2)h}.
\end{align*}
Notice that for $D_a=0$ this reduces to 
\begin{align*}
  \bc_{1-j}=-\bc_j+2\bc_{inj}.
\end{align*}
For $z=1$ and $k=m+j$, $j=1,\dots,p$, taking into account
\eqref{eq:bc2}, this extrapolation yields the value
\begin{align*}5\label{eq:bc22}
  \bc_{m+j}= \bc_{m+1-j},\quad  j=1,\dots,p.
\end{align*}

Fully discrete, high order, schemes are obtained by using an
appropriate Runge-Kutta ODE solver on \eqref{eq:weno-sd}.  To obtain
a fully explicit second order  scheme for approximations
$w^{n}_{i,j}\approx w_{i,j}(t_n)\approx w_i(x_j, t_n)$,
  $i=1,\dots,N$, $j=1,\dots,m$, we may use 
  \begin{equation} \label{eq:RK2-exp_old}
    \begin{aligned}
\boldsymbol{w}^{n+1/2}&=\boldsymbol{w}^n+\frac{\Delta t}{2}
\left(\boldsymbol{{\cal L}}^n+\boldsymbol{{\mD}}^{n} \right) \\
\boldsymbol{w}^{n+1}&=\boldsymbol{{w}}^n+\Delta t
\left(\boldsymbol{{\cal L}}^{n+1/2}+\boldsymbol{{\mD}}^{n+1/2}
\right).
\end{aligned}
\end{equation}
where 
we use the notation $\mD^*=\mD(\bw^*), \boc{L}^*=\boc{L}(\bw^*)$.
As observed in the previous section, the stability requirements of
this scheme would impose the use of time steps which are proportional
to the spatial step-size.

On the other hand, as shown in \cite{Donat1}, the stability
restrictions for 
 the following implicit-explicit Runge-Kutta 2 
(IMEX RK2) scheme
\begin{eqnarray}\label{imexrk21}
\boldsymbol{w}^{n+1/2}=\boldsymbol{w}^n+\frac{\Delta t}{2}
\left(\boldsymbol{{\cal L}}^n+\boldsymbol{{\mD}}^{n+1/2} \right) \\
\boldsymbol{w}^{n+1}=\boldsymbol{{w}}^n+\Delta t
\left(\boldsymbol{{\cal L}}^{n+1/2}+\boldsymbol{{\mD}}^{n+1/2}
\right). \label{imexrk22}
\end{eqnarray}
are the same as those of the purely hyperbolic case, i.e.
\begin{equation*}
\frac{u\Delta t}{\Delta z}\max_w \varrho (\bC'(w))\le C_1 \le 1.
\end{equation*}
Hence, for the ED model considered in this paper
\begin{equation}\label{IMEXstability}
\frac{\Delta t}{\Delta z}\le \frac{C_1}{u} 
\end{equation}
is sufficient for stability.

\subsection{IMEX-WENO schemes for the ED model}

The interplay between the variables $\bw$ and $\bc$ has to be taken
into account in order to solve the nonlinear system of equations
involved in an IMEX-WENO scheme such as
\eqref{imexrk21}-\eqref{imexrk22}. 
Specifically, for a single component chromatography,\eqref{imexrk21} explicitly reads as 
\begin{align}\label{eq:100}
\bw^{n+1/2}-\frac{\Delta t}{2}
\boldsymbol{C}^*(\bw^{n+1/2})\mathcal{A}
&=\bw^n+\frac{\Delta t}{2}
\boldsymbol{{\cal L}}(\bw^n).
\end{align}
This nonlinear equation should be solved by, e.g.,  Newton's
method. The difficulty that the nonlinearity in 
$\boldsymbol{C}^*(\bw^{n+1/2})$ is affected by the matrix
$\mathcal{A}$ can be overcome by performing a change of variables 
$\bc^{n+1/2}:=\boldsymbol{C}^*(\bw^{n+1/2})$, with which
\eqref{eq:100} reads now as:
\begin{align}\label{eq:100a}
\boldsymbol{W}^*(\bc^{n+1/2})-\frac{\Delta t}{2}
\bc^{n+1/2}\mathcal{A}
&=\mathcal{G}(\bw^n),\quad
\mathcal{G}(\bw^n):=
\bw^n+\frac{\Delta t}{2}
\boldsymbol{{\cal L}}(\bw^n),
\end{align}
where $\boldsymbol{W}^*(\bc)_{i,j}= W_i(c_j)$, $j=1,\dots,m, i=1,\dots,N$.

For simplicity in the description,
we shall consider first the case of a single component ($N=1$)
  and drop the index $n+1/2$. 
 In this case,  \eqref{eq:100a} becomes
\[ \left (1+ \frac{\eta}{1+bc_j}\right ) c_j - \theta
\left ( c_{j+1}-2 c_j + c_{j-1} \right )
- \mathcal{G}(\bw^n)_j =0
\] with (cf. \eqref{eq:207})
\begin{equation*}  \theta=\frac{\Delta t}{2}\mu=\frac{D_a \Delta t}{2 \Delta z^2}.
\end{equation*}
Thus, 
 the  vector of unknowns $\bc=(c_1,\ldots,c_m)^T$ satisfies the
 nonlinear system
\begin{equation}\label{eq:201}
  \mathcal{F}(\bc):=M(\bc)\bc-\mathcal{G}(\bw^n)=0,
\end{equation}  
where  the
matrix $M(\bc)$ is defined as
\[M(\bc)=E(\bc)-\mathcal{A}\]
with $E(\bc)$  diagonal 
\[ E_{jj}(\bc)= e(c_j)=1+
\frac{\eta}{1+bc_j}.  \]
We may solve this nonlinear system by  the standard  
Newton's method
\begin{align}\label{eq:206}
  \bc^{(\nu+1)}=\bc^{(\nu)}-\mathcal{F}'(\bc^{(\nu)})^{-1}\mathcal{F}(\bc^{(\nu)}),
  \nu=0,\dots, \quad \bc^{(0)}=\bc^{n},
\end{align}

Since
$\mathcal{F}_i(c)=e(c_i)c_i-\sum_{k=1}^{m} \mathcal{A}_{i,k}c_k$ and
$(e(c)c)'=1+\eta(1+bc)^{-2}$
\[\frac{\partial
  \mathcal{F}_i}{\partial
c_j}=\delta_{i,j}(e(c_i)c_i)'-\mathcal{A}_{i,j} \Rightarrow \mathcal{F}'(\bc)=\hat E (\bc)-\mathcal{A}, \quad \hat E= \text{diag}(1+ \frac{\eta}{(1+bc_i)^2}).
\]
Hence, solving
\eqref{eq:201} for each new
iterate only involves the solution of a tridiagonal system.

For the $N$-component ED model, with a Langmuir type adsorption
isotherm, we have
\begin{equation}\label{eq:204}
w_{i, j}=c_{i, j}
\left(1+\frac{\eta_i}{p_j}\right), \quad
p_j=1+\bb^T \bc_j
\end{equation}
where $i=1, \dots , N$ refers to the component of the mixture and $j=1,
\dots , m$ refers to the grid point under consideration.
  For an $N\times m$ matrix $A$, we denote $\mathcal{V}(A)\in
  \mathbb R^{Nm}$ given by juxtaposition of columns, i.e.,
  $\mathcal{V}(A)_{N(j-1)+i}=A_{i,j}$. By applying
  $\mathcal{V}$ to \eqref{eq:100a} and using the identity
  $\mathcal{V}(B X A^T)= (A\otimes B)\mathcal{V}(X)$, we get the vectorial equation

 \begin{equation}\label{matrixIMEX} 
{\mathcal{F}(\clC):=\mathcal{M}}(\clC) \,  {\clC}-{\mathcal{V}(\mathcal{G}(\bw^{n}))}=0
\end{equation}
for the unknown $\clC:=\mathcal{V}(\bc)$, where    ${\cal
  M}(\clC)$ is the $N m \times N m$ 
block tridiagonal matrix: 
 \begin{equation} 
{\cal M} (\clC)= {\cal E}(\clC) - \mathcal{A}\otimes I,
\end{equation}
with the $N\times N$ identity matrix $I$ and ${\cal E}(\clC)$ being a block-diagonal matrix 
\[E(\clC)= \begin{bmatrix} E^1 & {0}&{0}&{\dots}&{0}\\ 
{0}&{E^2}&{0}&{\dots}&{0} \\ 
{\dots}&{\dots}&{\dots}&{\dots}&{\dots} \\ 
{0}&{\dots}&{0}&{E^{m-1}}&{0} \\ 
{0}&{\dots}&{0}&{0}&{E^m}  \\ 
\end{bmatrix}, \,
\mathcal{A}\otimes I=\theta\begin{bmatrix} 
{-2 I}&{ I}&{0}&{\dots}&{0}\\ 
{ I}&{-2 I}&{ I}&{\dots}&{0} \\ 
{\dots}&{\dots}&{\dots}&{\dots}&{\dots} \\ 
{0}&{\dots}&{ I}&{-2 I}&{ I} \\ 
{0}&{\dots}&{0}&{ I}&{-2 I}  \\ 
\end{bmatrix}\]
with  $E^k$ ($k=1,\cdots,m$)  $N\times N$ diagonal matrices with
diagonal elements: 
\begin{eqnarray*}{E^k}_{ii}=1+\frac{\eta_i}{p_k} .
\end{eqnarray*}

As before, Newton's method requires  the computation of 
  $\mathcal{F}'(\boc{C})$ in \eqref{matrixIMEX}. From \eqref{eq:204}
  \begin{align*}
    \frac{\partial w_{i,k}}{\partial c_{j,l}}=\delta_{k,l}
    \left(
      \delta _{i, j}
\left(1+\frac{\eta_i}{p_k}\right)
-
    c_{i, k} b_j
\frac{\eta_i}{p_k^2}\right),\quad i, j=1, \dots , N, k, l=1, \dots , m.
\end{align*}
Therefore
$\mathcal{F}'(\boc{C})=\widehat{\mathcal{E}}-\mathcal{A}\otimes I$,
where $\widehat{\mathcal{E}}$ is the block diagonal matrix, whose
diagonal blocs 
are the  $N\times N$ (full) matrices
\begin{equation*}
(\hat
E^k)_{ij}=\left(1+\frac{\eta_i}{p_k}\right)\delta_{ij}-\frac{\eta_i
  b_j}{p^2_k} c_{i, k}, \quad i, j=1, \dots , N, k=1,\dots,m.
\end{equation*}

Thus, in the multi-component case, finding $\clC$ in
\eqref{matrixIMEX} by Newton's method \eqref{eq:206}
involves solving a 
block-tridiagonal system with small
blocks of size $N \times N$ at each iteration step. This can be
efficiently carried out by using a standard block tridiagonal LU
factorization algorithm (see \cite{Golub} for details).

The second step in the  IMEX-RK2 scheme \eqref{imexrk21} is
explicit, so that we directly obtain ${\boldsymbol w}_j^{n+1}$,
$j=1,\ldots,m$.  Through the computation of the only positive root
of  $R_{\bw}(p)=0$ with $\bw=\bw_j$ (e.g. by Newton's
method) we  obtain $p_j$ from which we easily get
the vector $\boldsymbol{c}_j^{n+1}$ at the next time step via
\eqref{eq:cw-def}.


Applying IMEX-WENO schemes to the ED model can, thus, be carried out
with a moderate computational effort. For multi-component mixtures,
the larger time steps allowed by the stability restrictions compensate
the additional effort with respect to a fully explicit alternative.

 
\section{Numerical experiments} \label{sec:numex}
In this section we perform several numerical experiments to
illustrate the behavior of the proposed IMEX-WENO scheme.

First, we consider a simple experiment that shows that the IMEX
alternative is able to obtain robust and reliable results when $D_a>0$
with a stability restriction of the type stated in
\eqref{IMEXstability}. Second, we consider the simulation of a
three-component mixture proposed in \cite{Javeed}.

%

\subsection{IMEX  versus Explicit Stability Restrictions}

We consider the same test for the single-component ED model considered
in section \ref{subsec:single} ($u=a=b=\eta=1, \epsilon=0.5$ and the
same  initial conditions).

\begin{figure}[ht]   
\begin{center} 
\includegraphics[scale=0.25]{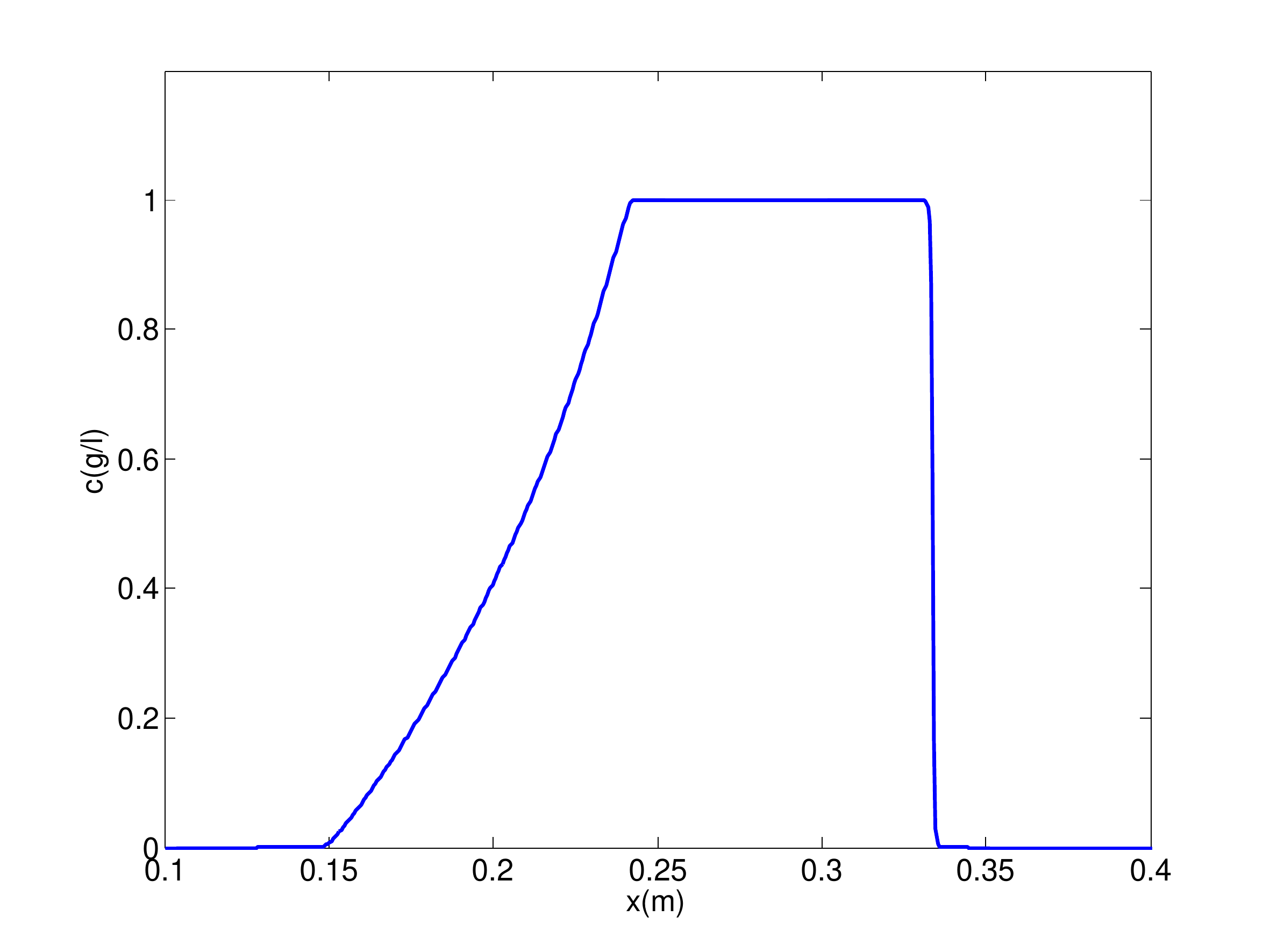}
\includegraphics[scale=0.25]{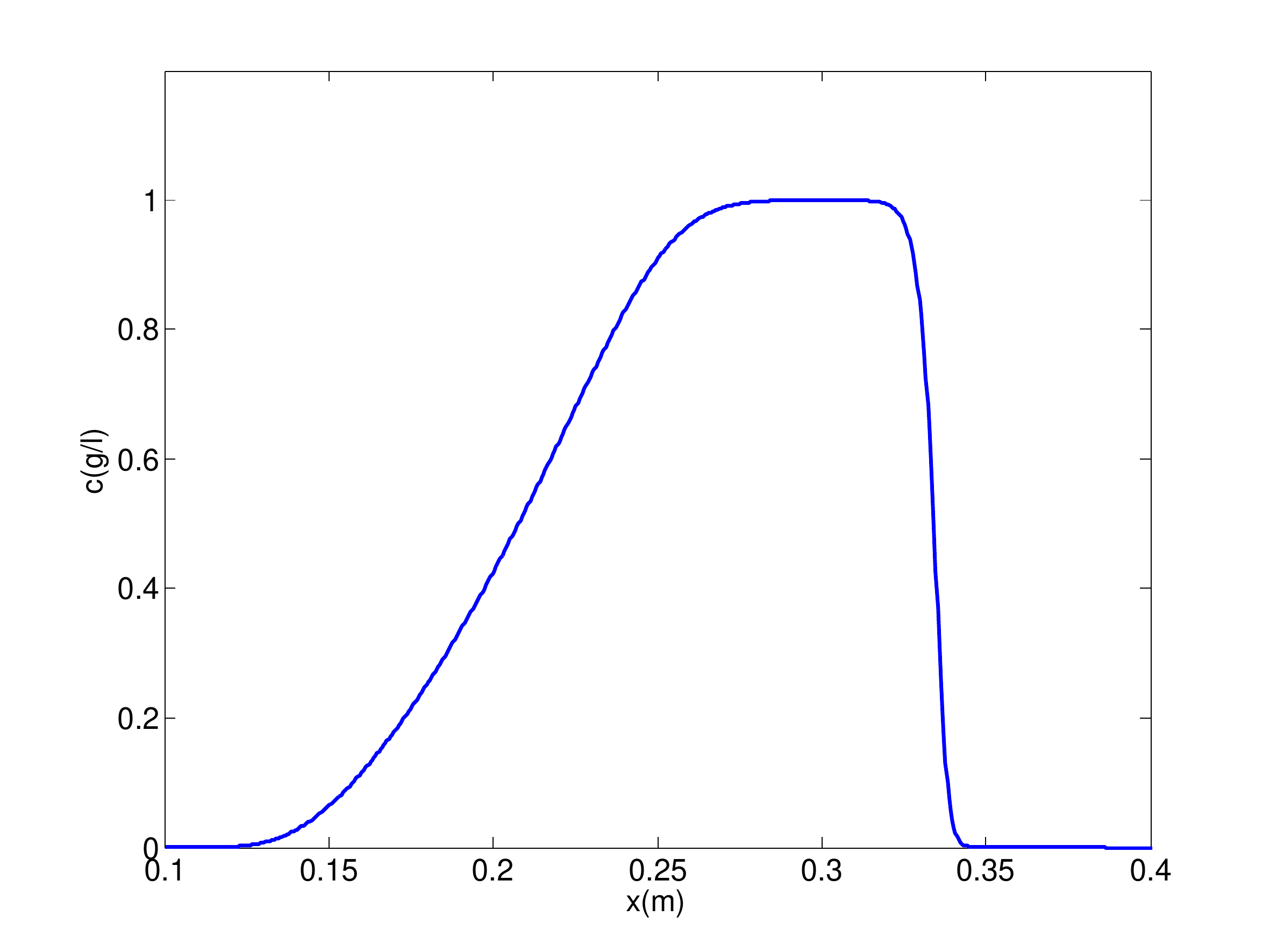}
\end{center}   
\caption{Single-component elution. $m=2000$,  $T=0.5$,
  $\Delta t/\Delta z=0.9$.  Left:  $D_a=0$ (Explicit). Right: $D_a=0.005$, IMEX-RK2 }   
\label{stability_sin}  
\end{figure}  

 In Figure
\ref{stability_sin} (left)  we show the numerical solution of the purely
hyperbolic problem at $T=0.5$ for $\Delta t/\Delta x =0.9$.
In Figure \ref{stability_sin} (right) we show the solution of the ED
model for $D_a=0.005$ obtained with the second order IMEX-WENO
scheme described in the previous section, also for  $\Delta
t/\Delta z =0.9$ and $T=0.5$. As expected, the 
simulation  is free of numerical artifacts and describes correctly the
effect of the parabolic term with respect to the solution obtained for
the purely hyperbolic model.
 
\begin{figure}[ht]   
\begin{center} 
\includegraphics[scale=0.25]{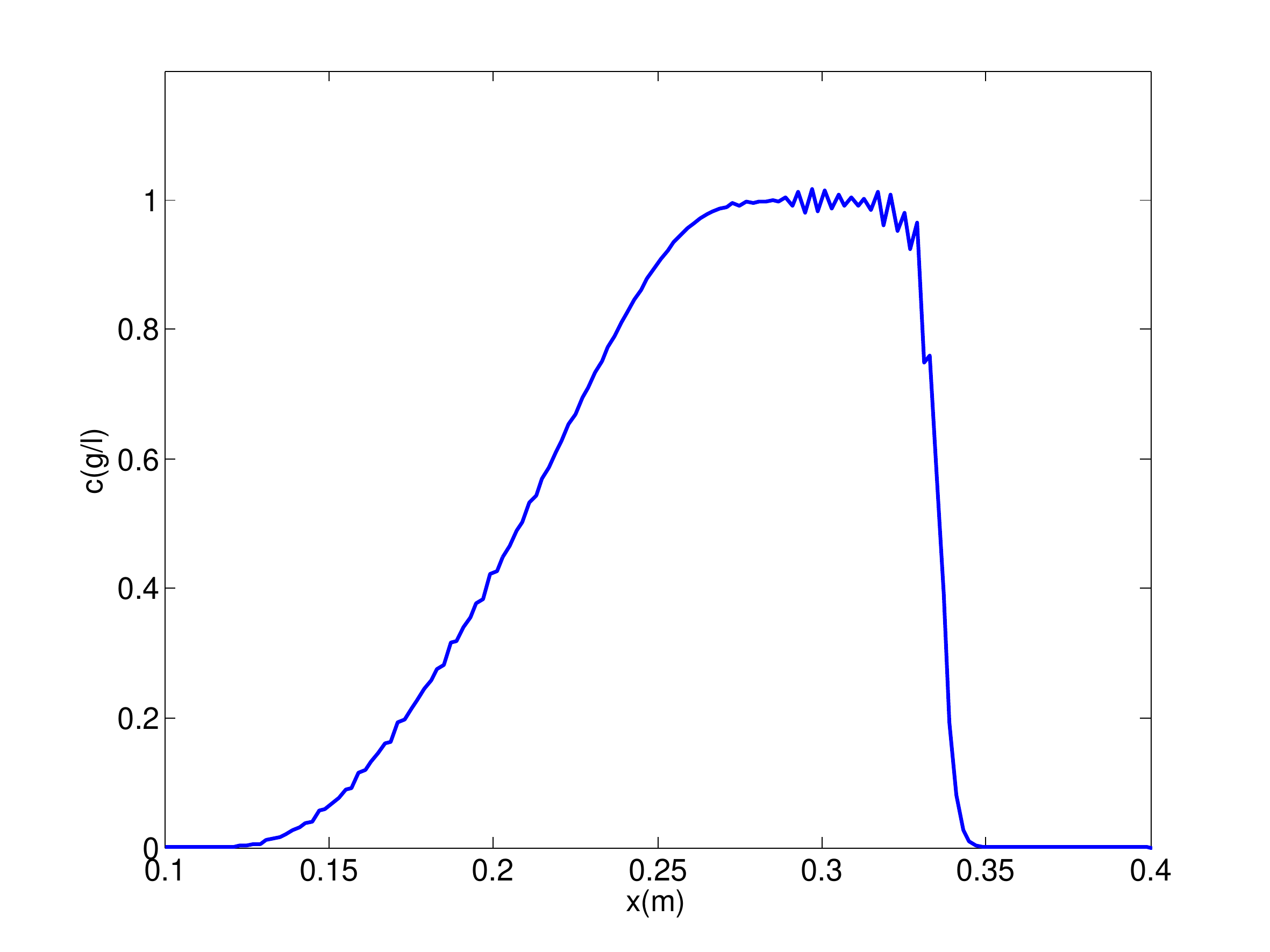}
\includegraphics[scale=0.25]{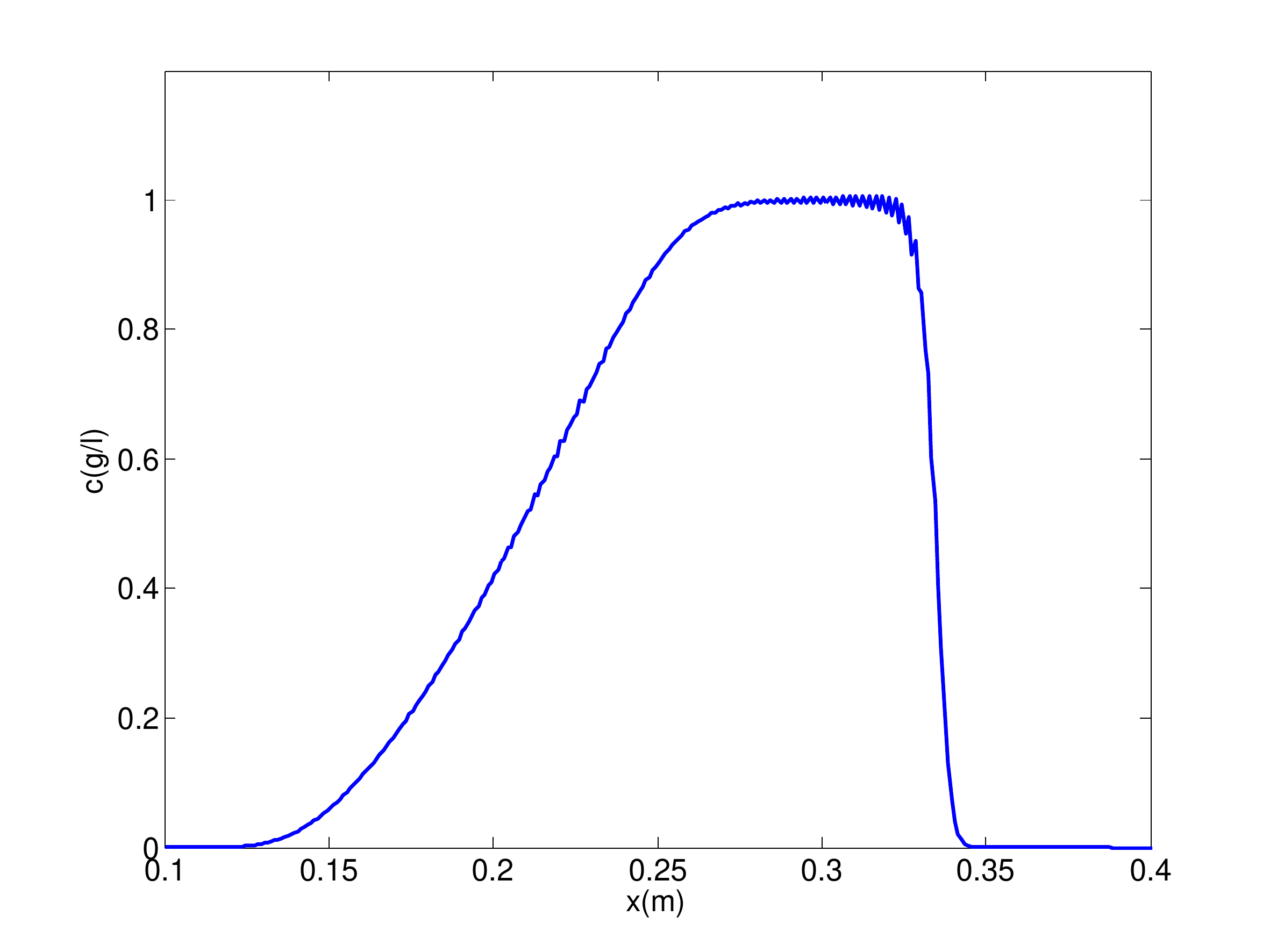}
\end{center}   
\caption{Single-component elution. Explicit
  schemes.   $D_a=0.005$, $T=0.5$. Left:  $m=500$, 
  $\Delta t/\Delta z=0.9$. Right:  $m=1000$, 
 $\Delta t/\Delta z=0.7$  }   
\label{stability_con}  
\end{figure}

As observed in section \ref{subsec:single}, the use of fully explicit
schemes when $D_a>0$ lead to     a 
stability restriction that depends on $\Delta z$.
In Figure \ref{stability_con} we show the numerical results obtained 
with the explicit two-step Runge-Kutta scheme \eqref{eq:RK2-exp_old}
for $m=500$ ( $\Delta t/\Delta x =0.9$) and $m=1000$ 
($\Delta t/\Delta x =0.7$). In both cases, the
stability restriction  \eqref{conditionstability} does not hold,
 and numerical oscillations develop. For $m=1000$ and  $\Delta
 t/\Delta x =0.9$ the numerical oscillations are so large that the
 numerical solution (not shown) is not representative of any
 meaningful behavior. 

In  Table \ref{tabla} we display the maximum values  for the ratio
$k/h$ for which no oscillatory behavior is observed in the numerical
simulation obtained with the fully explicit scheme \eqref{eq:RK2-exp_old}
(less than 3\% 
with respect to value of a reference solution, obtained with the IMEX
scheme for $m=2000$), and the value of the denominator in
\eqref{conditionstability}. The table shows that  condition
\eqref{conditionstability}, which was 
developed in \cite{Donat1} for the upwind first order scheme,  must be
enforced   in order to get stable numerical solutions when using any
explicit scheme for the ED model.

\begin{table}[htbp]
\begin{center}
\begin{tabular}{|c|c|c|c|c|}
\hline
m & 100 & 500 & 1000 & 2000 \\ \hline
 $(\Delta t/\Delta x)_{max}$ & 0.9 & 0.7 & 0.5 & 0.4 \\ \hline
$u+2 D_a/\Delta z$ & 1.1  & 1.5 & 2 & 3 \\ \hline 
$\Delta t/(u/\Delta z+2 D_a/\Delta z^2)$ & 0.99  & 1.05 & 1 & 1.2 \\ \hline 
\end{tabular}
\caption{Observed maximum values of $k/h$ for different number of cells $m$.}
\label{tabla}
\end{center}
\end{table}

\subsection{Three component displacement chromatography}

Displacement chromatography relies on the idea that one 
component  (the displacer) has  a
stronger affinity to the solid phase than any of the other components
in the sample mixture, hence it has the capability to displace the other
components of the mixture from the stationary phase.  For a
sufficiently long column and appropriate 
adsorption isotherms, the concentrations of the components form
rectangular regions of high concentration of one component in the
mixture. The series of such zones  are the
so-called isotachic train \cite{cazes2001encyclopedia}.

We consider  the case of a mixture of two components and one
displacer proposed in \cite{Javeed} (section
4.3). The values of the parameters are:
$a_1=4, a_2=5, a_3=6, b_1=4,b_2=5,b_3=1$. In addition, $N_t=10000$,
$\epsilon=0.5$ and $u=0.2$. 

 To compute the temporal
evolution of the 
concentrations of the two components, $c_1$ and $c_2$, and the
displacer, $c_3$  as they move along the column, we use the proposed
WENO-IMEX-RK2 scheme described in the previous section. \\[10pt]

{\em Experiment 1}.  Components 1 and 2 are injected between $t=0$ and
$t=0.1$ with $c_1=c_2=1$ at $z=0$. Component 3, the displacer, is
injected from $t=0.1$ with $c_3=1$. 

\begin{figure}[ht]   
\begin{center} 
\includegraphics[scale=0.250]{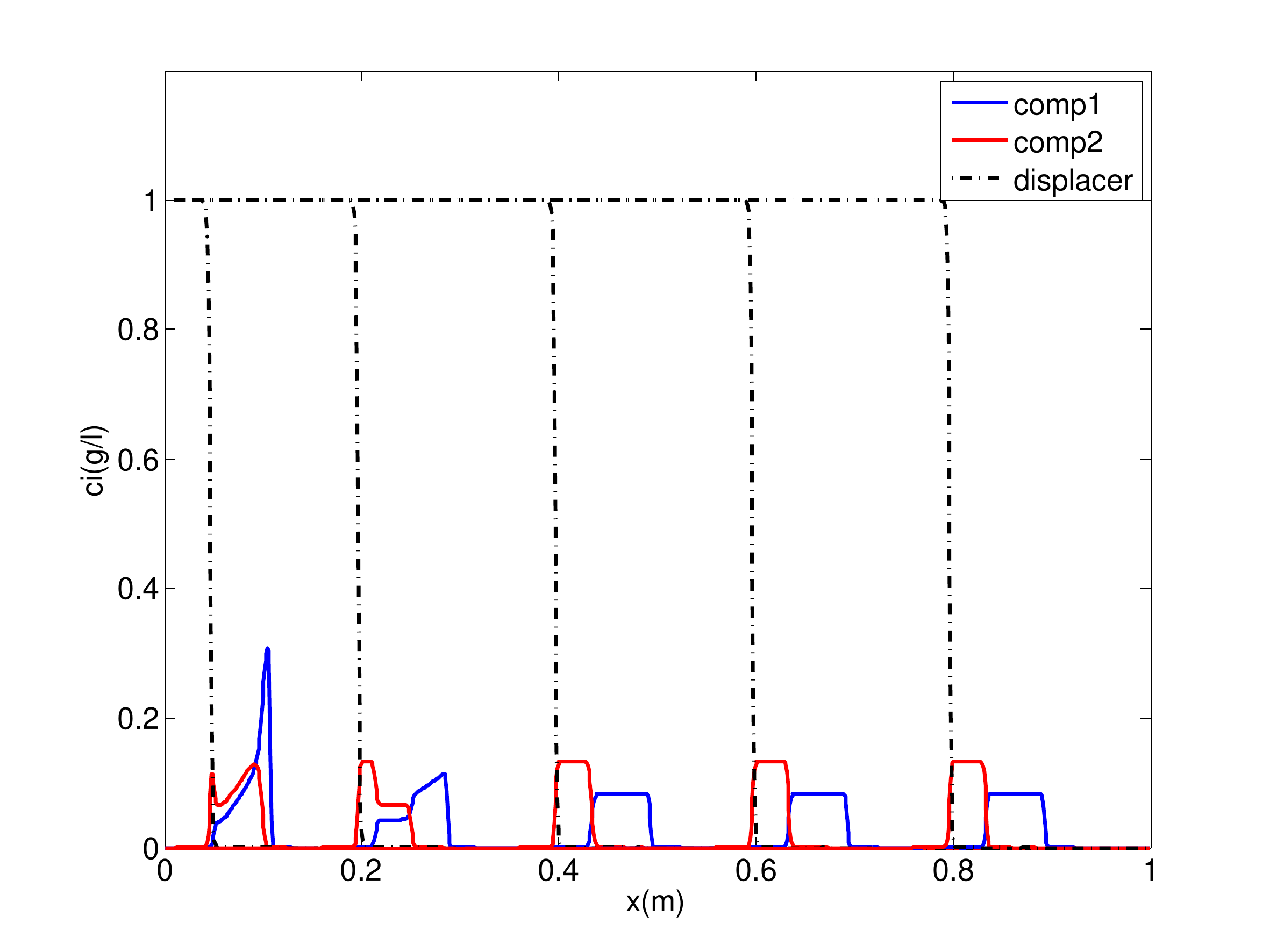}
\includegraphics[scale=0.250]{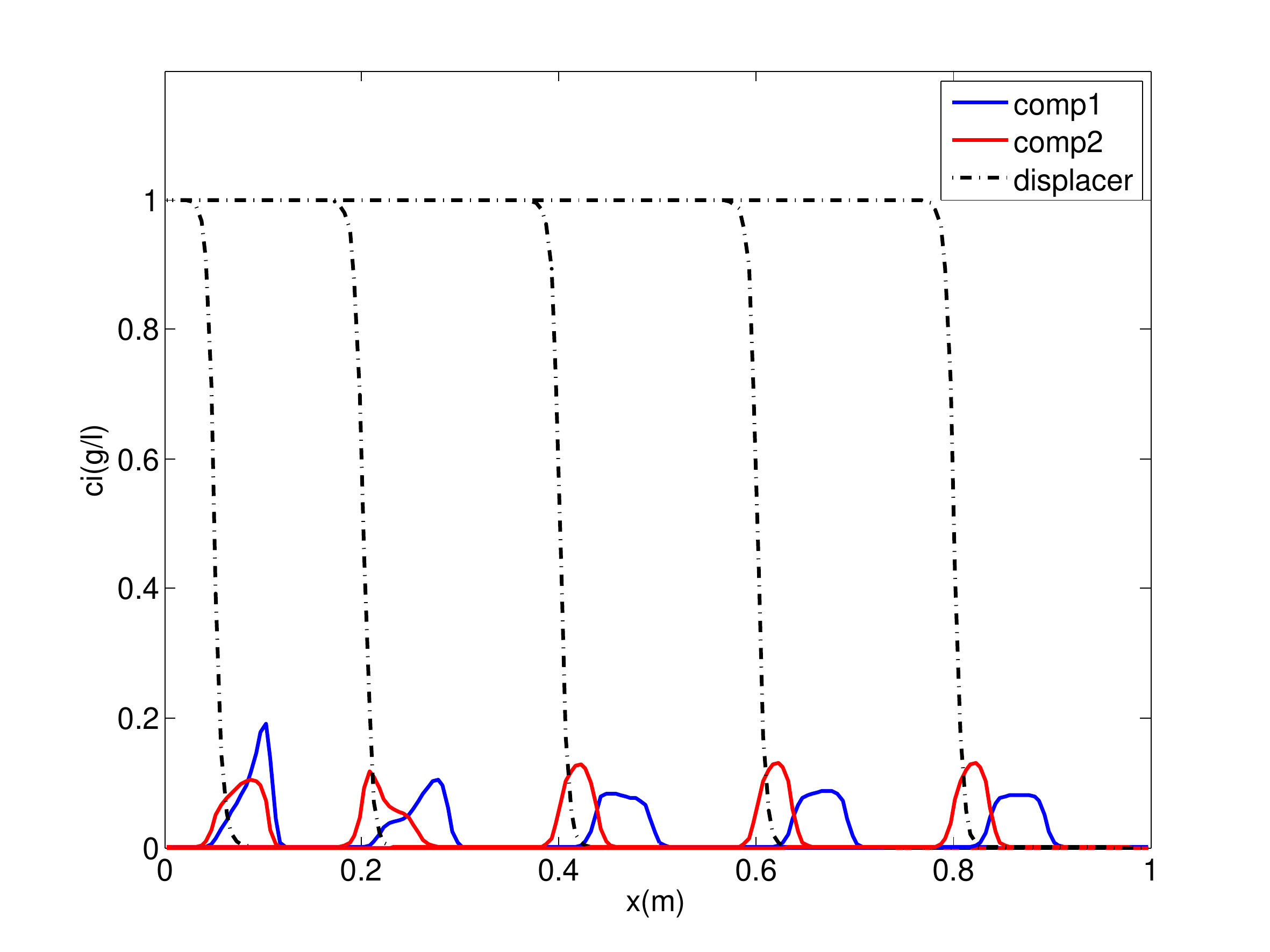}
\end{center}   
\caption{3-component test. Experiment 1: Numerical solution obtained
  with IMEX-RK2 $ \Delta t/\Delta z=4.0$.  Components 1, 2 and 3
  are shown in blue, red and black, respectively. Times: $T=1$, $T=4$,
  $T=8$, $T=12$ and $T=16$. Left $m=1000$. Right $m=200$.}   
\label{comp3disperso1}  
\end{figure}

 Figure \ref{comp3disperso1} shows the results for the first
 experiment for a CFL= $u \Delta t/\Delta z= 0.8$ and $m=1000$
 (left plot)
or  $m=200$ (right plot).  The
 formation of the displacement train can be clearly appreciated in
 both simulations. The results of the simulation for $m=1000$ can be
 compared with those reported in \cite{Javeed}. \\[10pt]

\begin{figure}[ht]   
\begin{center} 
\includegraphics[scale=0.250]{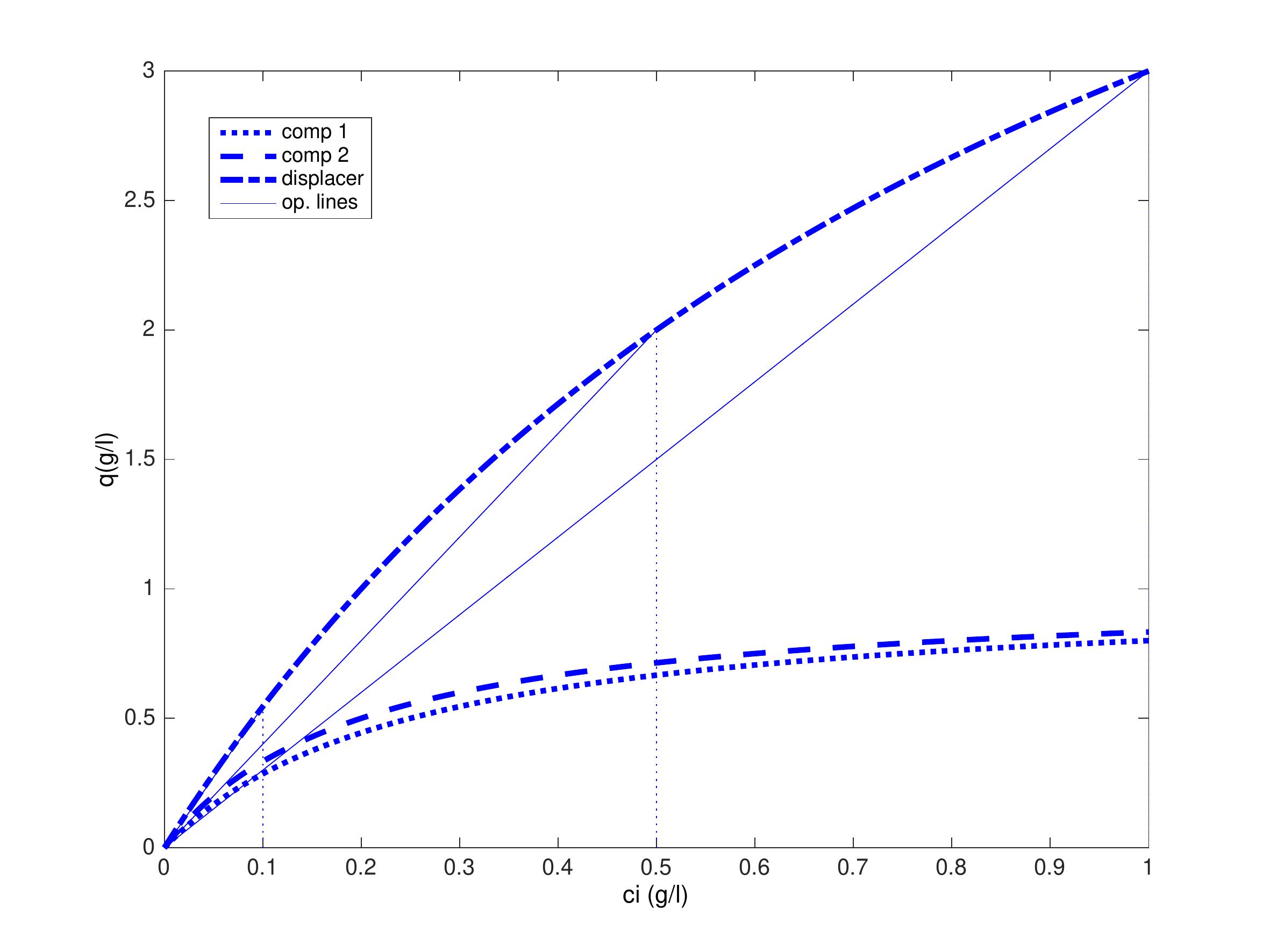}
\includegraphics[scale=0.250]{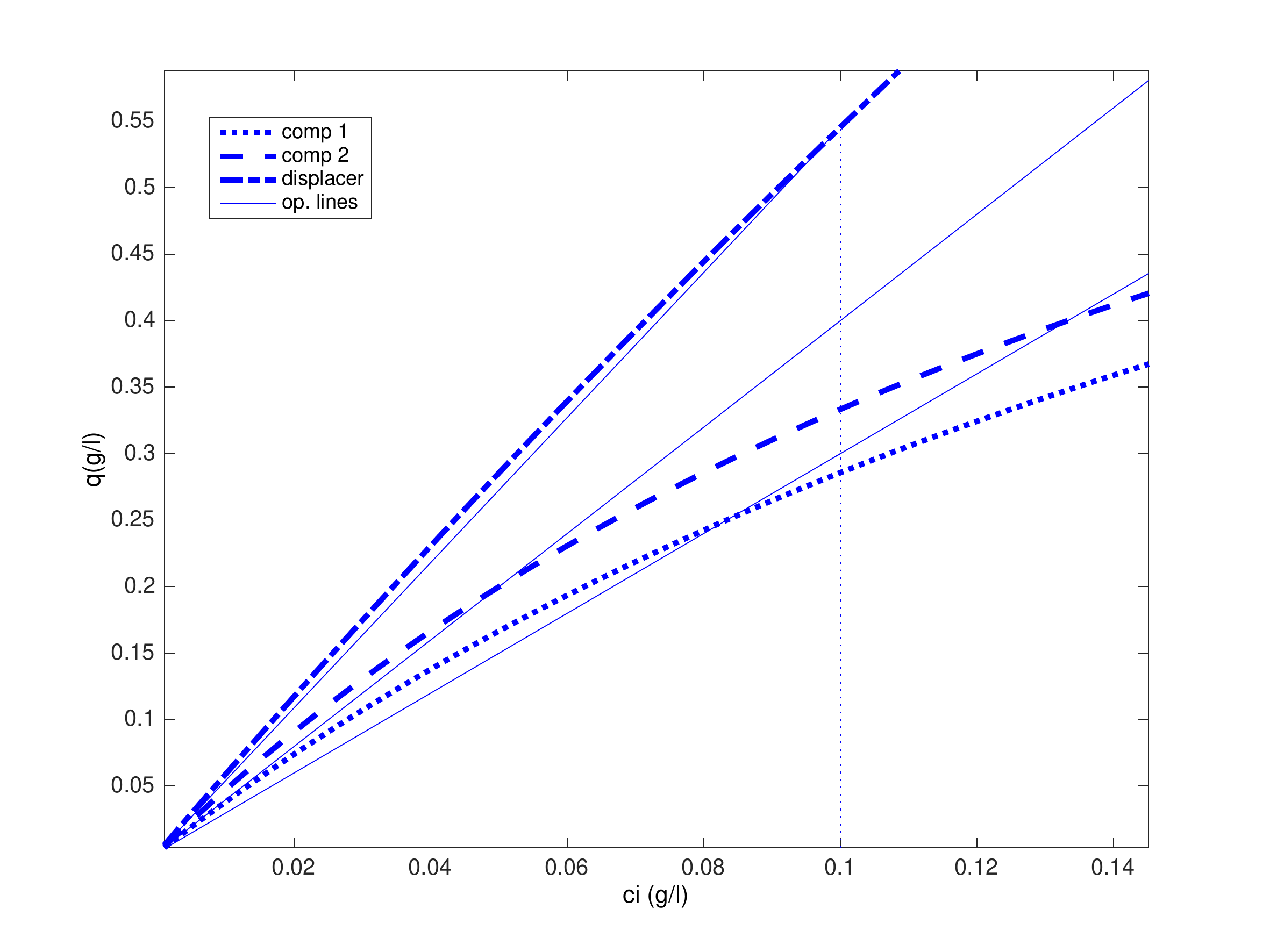}
\end{center}   
\caption{Left:Operating lines and isotherms . Right: Zoom of the area
  with the intersections.}   
\label{operatinglines}  
\end{figure}

\begin{figure}[ht]   
\begin{center} 
\includegraphics[scale=0.250]{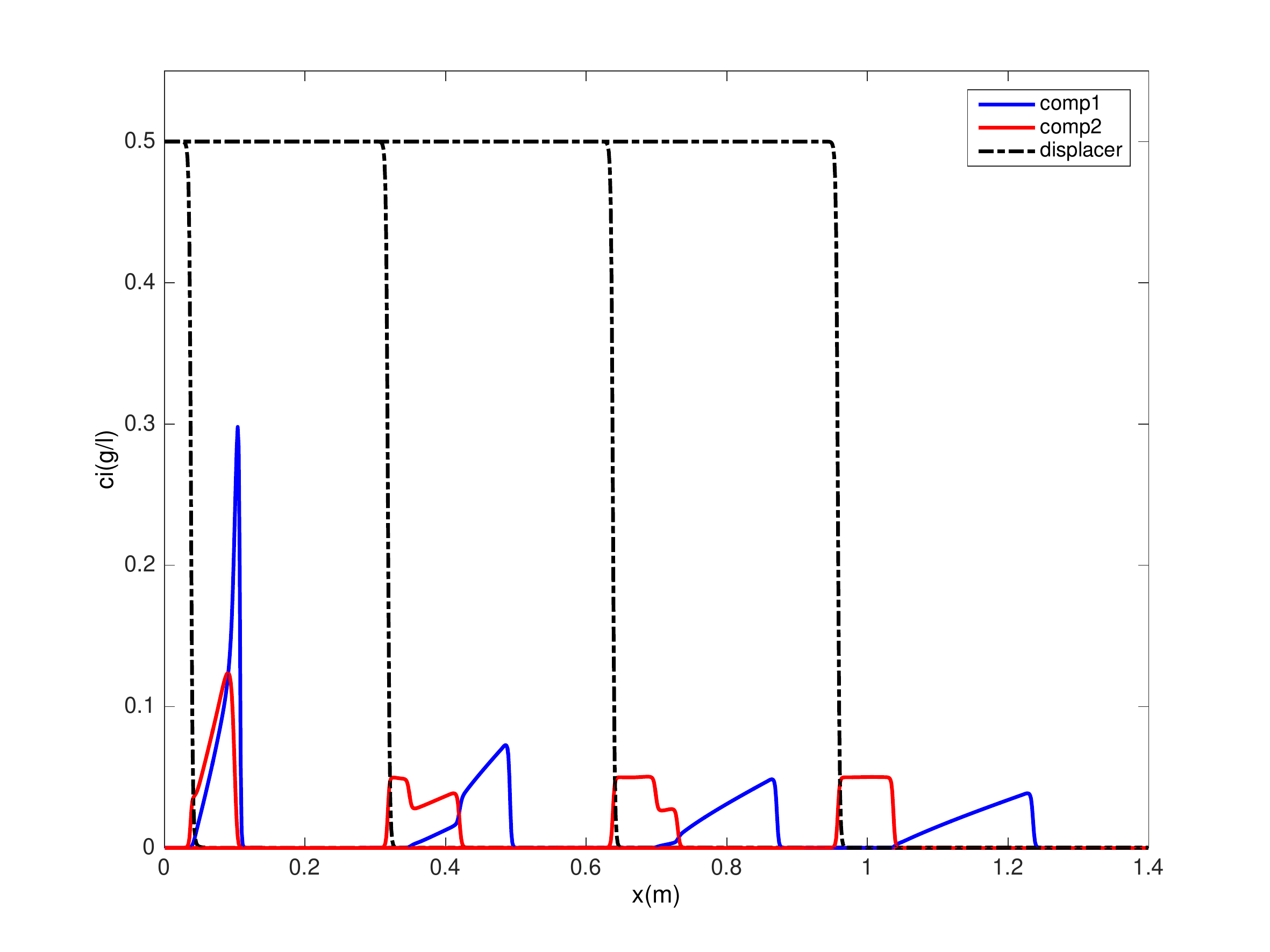}
\includegraphics[scale=0.250]{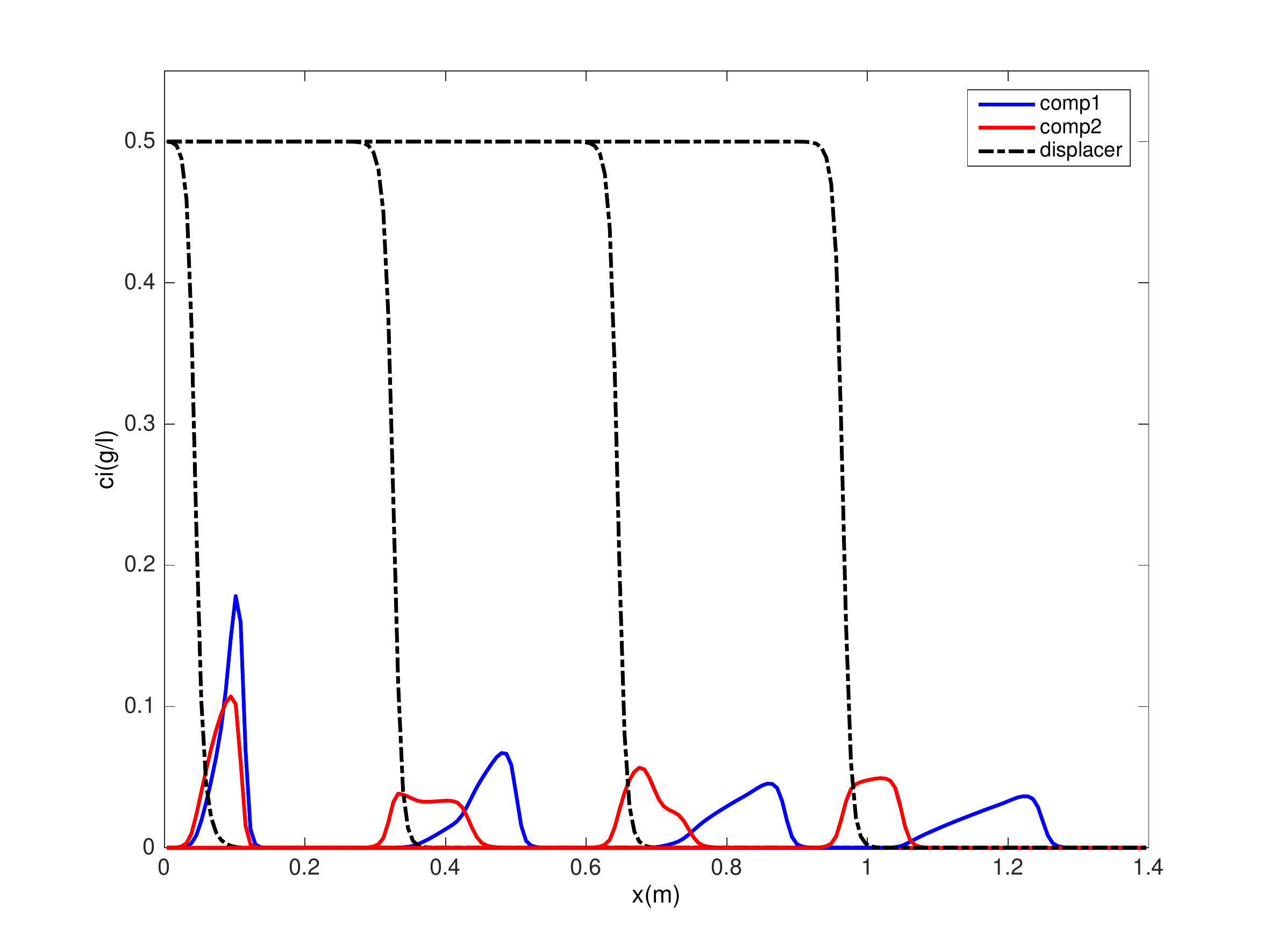}
\end{center}   
\caption{3-component test. Experiment 2: Numerical solution obtained
  with IMEX-RK2 $ \Delta t/\Delta z=4.0$.  Components 1, 2 and 3
  are shown in blue, red and black, respectively. Times: $T=1$, $T=8$,
  $T=16$ and $T=24$. Left $m=1000$. Right $m=200$.}   
\label{comp3disperso2}  
\end{figure}

{\em Experiment 2}: The concentration injected for the displacer is
$c_3=0.5$ and the values of the rest of the parameters are the same as
in experiment 1. This fact (see \cite{Javeed}) prevents the formation
of a rectangular pulse for component 1. As we can see in figure
\ref{operatinglines}, the line from the origin intersects the isotherm
of the second component, but not the corresponding to the first
component. This makes that only the second component can form a
rectangular pulse. Figure \ref{comp3disperso2} shows the time
evolution obtained with the IMEX-RK2 scheme. We see that the numerical
solution behaves as expected. \\[10pt]

\begin{figure}[ht]   
\begin{center} 
\includegraphics[scale=0.250]{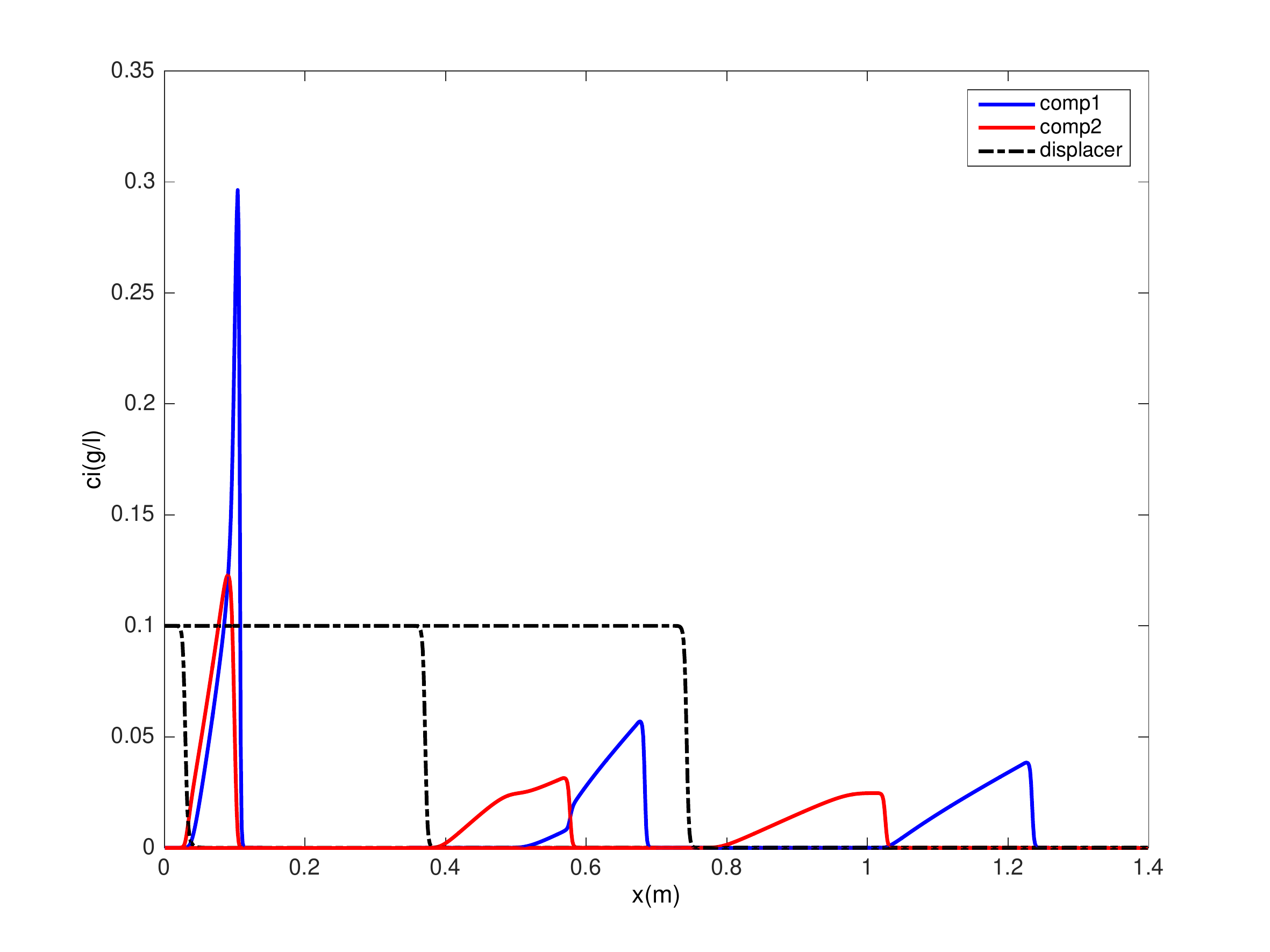}
\includegraphics[scale=0.250]{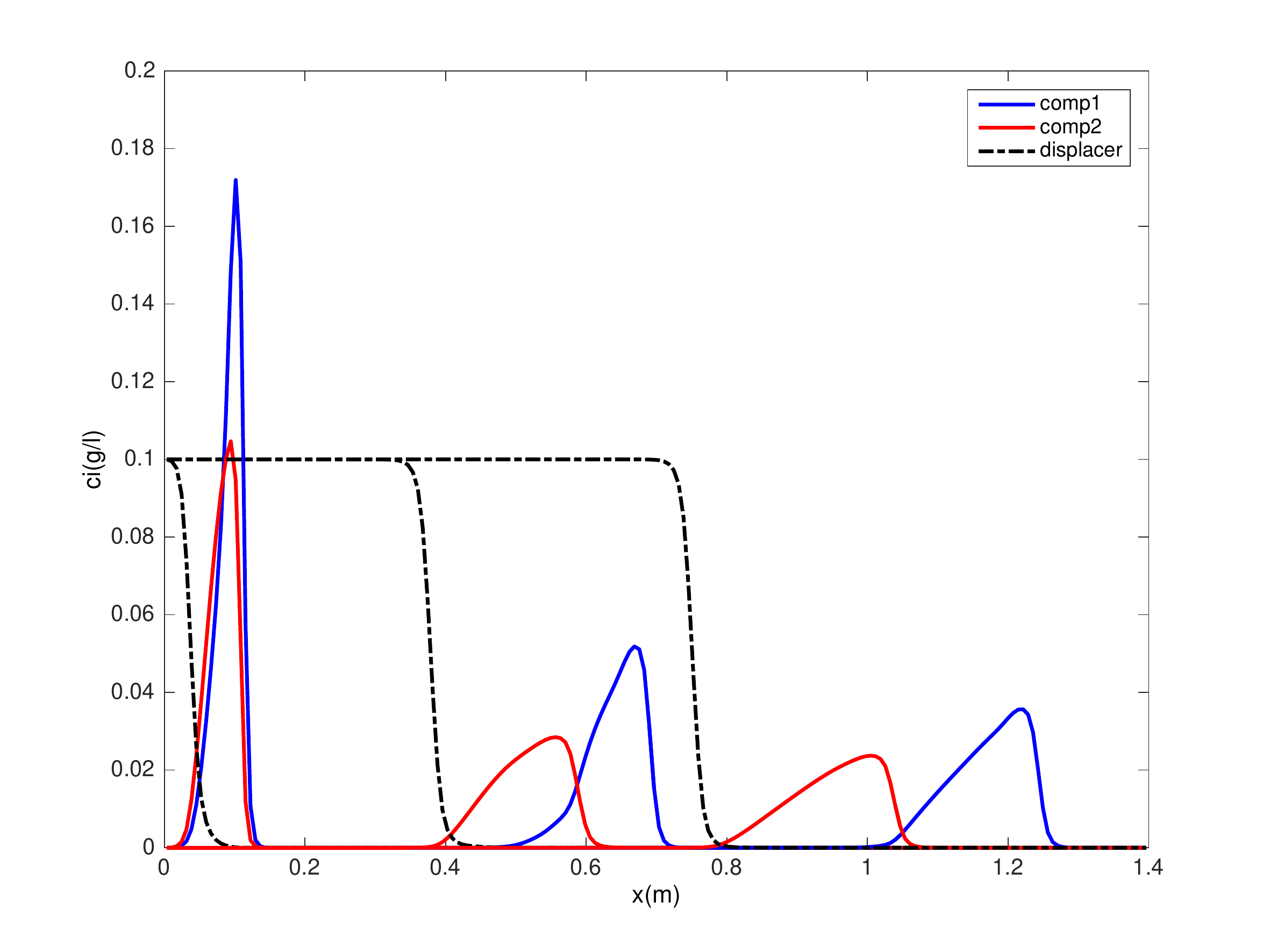}
\end{center}   
\caption{3-component test. Experiment 3: Numerical solution obtained
  with IMEX-RK2 $ \Delta t/\Delta z=4.0$.  Components 1, 2 and 3
  are shown in blue, red and black, respectively. Times: $T=1$, $T=12$ and $T=24$. Left $m=1000$. Right $m=200$.}   
\label{comp3disperso3}  
\end{figure}

{\em Experiment 3}: The concentration injected for the displacer is
further reduced, to the
value $c_3=0.1$. The values of the rest of the parameters are the same
as in experiment 1. 

In this case, as can be seen in Figure \ref{operatinglines}, none of the isotherms is intersected by the operating line. According to this, both components fail to form equilibrated rectangular pulses. 
The results are shown in Figure \ref{comp3disperso2}. Again, our IMEX-RK2 scheme reproduces correctly the expected behavior.

\section{Conclusions} \label{sec:conc}

In this paper we have examined the ED model with Langmuir-type
adsorption isotherms. We have proven that for $D_a=0$, the model can be
written as a system of conservation laws, since  the relation between
the conserved 
variables $\bw= \bc+\frac{1-\epsilon}{\epsilon}\bq(\bc)$ and the
physical concentrations, $\bc$, admits a smooth, globally well defined
inverse, $\bC(\bw)$. The properties of the function
$\bW(\bc)=\bc+\frac{1-\epsilon}{\epsilon}\bq(\bc)$,  relating the
conserved variables 
and the physical concentrations for Langmuir-type isotherms, are
exploited in order to show that the system of conservation laws is
strictly hyperbolic. 

The inverse function $\bC(\bw)$ does not admit an explicit expression
for $N>1$, however we show that $\bC(\bw)$ can be efficiently computed
by finding the only  
positive root of a  rational function.
The capability to compute $\bC(\bw)$ at a reasonable cost allows us to
design fully conservative explicit schemes, which ensure the correct
propagation of discontinuous, or nearly discontinuous, fronts. We
discuss the advantages of the fully conservative schemes compared to the
non-conservative  high-resolution
schemes proposed in \cite{Javeed}.

Implicit-Explicit strategies are often considered  for
convection-dominated parabolic systems of PDEs, such as the ED model
for $D_a<<1$, since  the stability restrictions of these schemes are
less severe that those obtained for fully discrete schemes. We propose
 a second  order WENO-IMEX 
scheme in order to illustrate the numerical difficulties arising from the
non-explicit relation between the conserved and physical variables in
the ED model. We show that the structure of the Langmuir adsorption isotherms
can be used to set up a change of variables that allows us to compute
the physical concentrations (instead of the 
conserved variables) in implicit steps. The solution of the
resulting nonlinear system  of equations that arises at each time step
 can be found by Newton's method and it involves solving a
block-tridiagonal linear system, with small blocks that are of the size of the
number of components in the mixture, at each iteration. These systems
may be efficiently 
solved by standard block-tridiagonal routines, making the proposed
WENO-IMEX-RK technique a very efficient and robust scheme for this
model. Several numerical experiments show the performance and capabilities of
the proposed numerical scheme. 

 In this paper we have applied the scheme to the equilibrium dispersive model for a fixed bed. However,  it might be adapted to more general models of chromatographic column and also to simulating moving beds (SMB) models.

\section*{Acknowledgements}
This research was partially supported by Ministerio de Econom\'ia y Competitividad under grant  MTM2014-54388-P with the participation of FEDER.

 \end{document}